\RequirePackage{fix-cm}
\documentclass[smallextended]{svjour3}       
\smartqed  
\usepackage{graphicx}

\usepackage{placeins}

\usepackage[pagewise]{lineno}
\usepackage{lipsum,blindtext}
\usepackage{mathrsfs}
\usepackage{tabularx}
\usepackage{multirow}
\usepackage{tabu}
\usepackage{menukeys}
\usepackage{attachfile}
\usepackage{subcaption}
\usepackage{circuitikz}
\usepackage[ruled,vlined,linesnumbered]{algorithm2e}
\newcommand\numberthis{\addtocounter{equation}{1}\tag{\theequation}}
\usepackage{forest} 
\usepackage{fontawesome5}
\usetikzlibrary{arrows}
\usetikzlibrary{decorations.pathreplacing}
\usetikzlibrary{patterns}

\ctikzset{resistor = european}

\makeatletter 
\pgfdeclarepatternformonly[\LineSpace,\tikz@pattern@color]{my north east lines}{\pgfqpoint{-1pt}{-1pt}}{\pgfqpoint{\LineSpace}{\LineSpace}}{\pgfqpoint{\LineSpace}{\LineSpace}}%
{
    \pgfsetcolor{\tikz@pattern@color} 
    \pgfsetlinewidth{0.4pt}
    \pgfpathmoveto{\pgfqpoint{0pt}{0pt}}
    \pgfpathlineto{\pgfqpoint{\LineSpace + 0.1pt}{\LineSpace + 0.1pt}}
    \pgfusepath{stroke}
}
\makeatother 
\newdimen\LineSpace
\tikzset{
    line space/.code={\LineSpace=#1},
    line space=3pt
}

\usepackage{listings}
\usepackage{tikz}
\usetikzlibrary{arrows}
\usetikzlibrary{arrows,automata}
\usepackage{amsmath,amssymb}

\usepackage{float}

\usepackage[T1]{fontenc}
\usepackage{lmodern} 

\begin{document}

\title{Statistical static timing analysis via modern optimization lens
}
\subtitle{\textcolor{black}{I.~Histogram--based approach.}}

\titlerunning{SSTA via modern optimization lens\textcolor{black}{: I. Histogram--based approach}
}        

\author{Adam Bos\'ak         \and
        Dmytro Mishagli \and
        Jakub Mare\v{c}ek
}


\institute{A. Bos\' ak \at
Technical University of Denmark, Kongens Lyngby, Denmark\\
\and 
A. Bos\'ak and J. Mare\v{c}ek \at
              Czech Technical University in Prague, the Czech Republic\\          
           \and
           D. Mishagli \at
              University College Dublin, Belfield, Ireland \\
              \email{\{bosadam@seznam.cz, dmytro.mishagli@ucd.ie, jakub.marecek@fel.cvut.cz\}} 
}

\date{Received: date / Accepted: date}

\maketitle

\begin{abstract}
{\color{black}
Statistical static timing analysis (SSTA) is studied from the point of view of mathematical optimization. We present two formulations of the problem of finding the critical path delay distribution that were not known before: (i) a formulation of the SSTA problem using Binary--Integer Programming and (ii) a practical formulation using Geometric Programming. For simplicity, we use histogram approximation of the distributions. Scalability of the approaches is studied and possible generalizations are discussed.
}

\end{abstract}

\section{Introduction}

Integrated Circuits (ICs) must work at expected frequencies with respect to the timing constrains specified in their designs, which is checked by Computer--Aided Design (CAD) tools. The key challenge is to satisfy these constraints in an optimal way, so that the area taken by the designs on a chip and the power consumption are minimized. This leads to a particular class of optimization problems. The \textcolor{black}{development of the} techniques used for design verification forms a separate field of study, the Timing Analysis~\cite{sapatnekar2004Timing}.

With the decrease of features sizes, the impact of variations that occur in manufacturing processes, the process variations, increases. These variations lead to uncertainties in the parameters of the transistors and interconnects, which affect the delays and hence the overall performance of a circuit. For example, the performance of the same designs can differ from chip to chip (intradie or global variations). At the same time, a single design can have variation of delays in different parts of a die (inter--die or local variations). The most reliable way to take these variations into account in order to predict the yield (i.e., the fraction of correctly functional chips among all fabricated) is to run Monte Carlo (MC) simulations. For modern Very Large Scale Integration (VLSI) designs, this is very expensive, which again increases the cost of chips.

A less computationally expensive approach is to use a Static Timing Analysis (STA), which is still the most common way to take into account systematic (global) process variations~\cite{sapatnekar2004Timing,vygen2006SlackStaticTiming,bhasker2009StaticTimingAnalysis}. This approach treats variations as single and determined values (so-called \textit{corner} values). The delay values computed in such a way are too pessimistic~\cite{blaauw2008StatisticalTimingAnalysis}, which results in an increase in the cost of chips when these delays are mitigated~\cite{visweswariah2003DeathTaxesFailing}. An additional challenge is that variations have substantial non-Gaussian behaviour and are often strongly correlated. In modern ultra-VLSI circuits (5~nm and below), the impact of random correlated processes and fluctuations has become even more important. Thus, an alternative approach has been developed.

Statistical Static Timing Analysis (SSTA) addresses randomness in a natural way, treating delays in a system as Random Variables (RVs) from the very beginning. The analysis then allows us to determine the mean value of the delay across selected paths. The maximum delay corresponds to the critical path. Current industrial realizations of SSTA allow one to determine moments of delay distributions and/or their quantiles~\cite{chang2003StatisticalTimingAnalysis,chang2005StatisticalTimingAnalysis,chang2005ParameterizedBlockbasedStatistical}. In principle, SSTA can give a delay distribution of the whole circuit. This makes SSTA comparable to MC simulations in terms of accuracy. At the same time, SSTA algorithms are much less resourceful than MC but have higher computational complexity than deterministic STA.

This work is motivated by~\cite{freeley2018StatisticalSimulationsDelay,mishagli2020RadialBasisFunctions}, where an approach was proposed to deal with non-Gaussian distributions of gates' delays without loss of information while keeping complexity low enough. The authors based their method on the exact solution to the problem of finding a distribution of a logic gate delay assuming that all distributions are Gaussian. The resulting non-Gaussian distribution is then decomposed into a mixture of Radial Basis Functions with fixed shape parameters and locations, but unknown mixing coefficients (the weights). Such a mixture, which is referred to as \emph{Gaussian comb}, can be obtained by means of a Linear Programming. This is the main advantage of the proposed model.

\textcolor{black}{The aim of this paper is to investigate a classical SSTA problem from the mathematical optimization point of view, which is, to the best of our knowledge, done for the first time. For simplicity, we will represent actual distributions with histograms, which inevitably introduces accuracy drop. We make this choice deliberately, as our focus is on proofing the concept rather than developing an accurate approach. We give two formulations of the SSTA problem as an optimisation problem \textit{via} (i)~Binary--Integer Programming and (ii)~Geometric Programming. Finding efficient solutions will be the subject of a separate study where we combine the approach presented in this work with the Gaussian comb model.
}

The paper is organized as follows. In Section~\ref{sec:background}, \textcolor{black}{we give necessary mathematical preliminaries, terminology, and discuss related work. Section~\ref{sec:ssta_setup} discusses standard and straightforward (\textit{i.e.}, non-optimization) approaches to SSTA and gives the statement of the problem. We give our solutions \textit{via} optimization techniques in Sections~\ref{sec:ssta_bip} and~\ref{sec:ssta_gp}. Section~\ref{sec:ssta_bip} discusses formal mathematical solution by means of Binary Integer Programming.
Section~\ref{sec:ssta_gp} shows that the scalability can be improved by utilizing Geometric Programming. 
We then conclude the paper in Section~\ref{sec:discussion} with overall discussion of the obtained results.}

\section{Background and Related Work}\label{sec:background}

In this Section, we start from some mathematical preliminaries, followed by introducing the terminology, and finally briefly describe the key results in the Statistical Static Timing Analysis.

\subsection{Mathematical Preliminaries}
{
\color{black}
Here we will present only necessary definitions and
refer to Boyd~\textit{et al.}~\cite{boyd2007TutorialGeometricProgramming} for a very accessible overview of Geometric Programming. For further details, we refer the reader to the references to the literature cited there.}

Monomial is a function: $\mathbb{R}^n_{++} \rightarrow \mathbb{R}_{++} $\footnote{The domain of the monomials is the non-negative quadrant of $\mathbb{R}^n$. We assume that the optimum values cannot be zero, and therefore the domain is of the form $\mathbb{R}^n_{++}$.}

\begin{equation}
    f(x_{1},x_{2},...,x_{n}) = cx_{1}^{a_{1}...}x_{n}^{a_{n}},
\end{equation}  
where $c > 0$ and ${ a_i \in \mathbb{R} }$. We call $c$ the coefficient of~the~monomial and $a_{i}$ exponents of the monomial. We refer to a sum of monomials as a posynomial, that is, a function of the form

\begin{equation}
    f(x_{1},x_{2},...,x_{n}) = \sum_{k=1}^{K} c_{k}x_{1}^{a_{1k}...}x_{n}^{a_{nk}}.
\end{equation}

It is easy to see from the definitions that for the set of all monomials $A$ and for the set of all posynomials $B$, it holds ${ A \subseteq B }$. One should also note that posynomials are closed under addition, multiplication, positive scaling, and results in a posynomial, when divided by a monomial.

A Geometric Program (GP) is an optimization problem of the form
\begin{alignat*}{3}\label{GP}
 & \text{minimize} & f_0(\mathbf{x}) \\ \numberthis
 & \text{subject to} \quad & f_i(\mathbf{x}) & \leq 1, \quad & i &=1 ,..., n\\
                && g_i(\mathbf{x}) & = 1, \quad & i &=1 ,..., p,
\end{alignat*}
where $x_i$ are optimization variables, $f_i$ are posynomial functions, and $g_i$ are monomials. We call~\eqref{GP} a GP program in a standard form.

{\color{black}
An Integer Linear Programming (ILP) problem is written in general form as follows (see, \textit{e.g.}, Wolsey~\& Nemhauser~\cite{wolsey1999IntegerCombinatorialOptimization}):
\begin{alignat*}{3}\label{eq:MIP_general}
 & \text{maximize} & \mathbf{c}^T \mathbf{x}  \numberthis \\
 & \text{subject to} \quad &\mathbf{A} \mathbf{x} &\leq \mathbf{b}, \\
                && \mathbf{x} &\geq \mathbf{0}, \\
                && \mathbf{x} &\in \mathbb{Z}^n.\\  
\end{alignat*}
If the domain polyhedron is intersected with the hypercube $\{0, 1\}^n$, we talk about Binary--Integer Programming (BIP) problems. This can be specified using constraints


\begin{equation}\label{eq:BIP_constraints_general}
    \begin{aligned}
    \mathbf{0} \leq \mathbf{x} \leq \mathbf{1},\\
    \mathbf{x} \in \mathbb{Z}^n.\\
    \end{aligned}
\end{equation}
}

\subsection{Definitions}\label{sec:background:definitions}
Throughout the paper, we will use the following commonly accepted terminology~\cite{blaauw2008StatisticalTimingAnalysis}. A logic circuit can be represented as a timing graph $G(E,V)$, where the \textit{graph} and its \textit{paths} are defined as follows.

\begin{definition}
A timing graph $G(E,V)$ is an acyclic directed graph, where $E$ and $V$ are the sets of edges and vertices, respectively. The vertices correspond to logic gates of a circuit. The timing graph always has one source and one sink. The edges are characterised by weights $d_i$ that describe delays.
The timing graph is called a \emph{statistical timing graph} within SSTA when the edges of the graph are described by RVs.
\end{definition}

The task then is to determine the critical (longest) path.

\begin{definition}
Let $p_i (i = 1,\ldots,N)$ be a path of ordered edges from the source to the sink in a timing graph $G$ and let $D_i$ be the path length of $p_i$. Then $D_{\text{max}} = \max(D_1,\ldots,D_N)$ is called the \textit{SSTA problem of a circuit}.
\end{definition}

Since the gates have internal structure presented by corresponding combination of transistors, this results in a characteristic time needed for the gates to operate. This is one of the sources of delays in a circuit. Due to delays, input signals can have different arrival times, and therefore, the delay of a gate is determined by the maximum of input delays. On the other hand, the operation time of a gate can have a significant impact on the circuit delay, in addition to the arrival times. In such a case, the delay of a gate itself should be added to the result of the $\max$ function:
\begin{equation}\label{eq:gate_operation}
  d_{\text{gate}} = \max(d_1,d_2) + d_0 + d_{\text{int}} + \ldots,
\end{equation}
where $d_1$, $d_2$ are delays in input signals, $d_0$ is a gate delay (due to its operation time) and $d_{\text{int}}$ is an interconnect delay.

The calculation is straightforward in the case of a deterministic timing analysis, but is not the case when uncertainty arises. As we have already mentioned, within the SSTA, the arrival and gate operation times are described by RVs given by the corresponding distributions. Therefore, the delay \eqref{eq:gate_operation} can be written as
\begin{equation}\label{eq:gate_operation_rv}
    \zeta_{\text{gate}} = \max(\xi_1,\xi_2) + \xi_0 + \xi_{\text{int}} + \ldots,
\end{equation}
where $\xi_1$ and $\xi_2$ are RVs that describe the arrival times of input signals, $\xi_0$ and $\xi_{\text{int}}$ are RVs related to the gate operation time and the interconnect delay respectively. The whole gate delay, $d_{\text{gate}}$, is now an RV itself and is indicated by $\zeta_{\text{gate}}$. In principle, $\xi_0$ and $\xi_{\text{int}}$ can be combined in a single RV, thus, the latter will be omitted in future discussion.

Therefore, as we can see from \eqref{eq:gate_operation_rv}, two operations fully describe delay propagation at the gate level: \textcolor{black}{(i)~the maximum of delays entering a gate and (ii)~the summation of the latter with the delay of the gate. These operations are often called \textit{atomic operations} of SSTA (see, \textit{e.g.}, works by Cheng~\textit{et al.}~\cite{cheng2007NonLinearStatisticalStatic,cheng2012FourierSeriesApproximation}).} In the language of distributions, Eq.~\eqref{eq:gate_operation_rv} gives a \emph{convolution} of probability density functions of the RVs $\max(\xi_1,\xi_2)$ and $\xi_0$. In this work, we consider a histogram approximation to the problem, which will be discussed in the next sections.

\subsection{Related Work}

There were excellent reviews of the work done in the early stage of the SSTA era, 2001--2009, namely by Blaauw~\textit{et al.}~\cite{blaauw2008StatisticalTimingAnalysis} and Forzan~\textit{et al.}~\cite{forzan2009StatisticalStaticTiming}. A good overview is also conducted by Beece~\textit{et al.}~\cite{beece2014TransistorSizingCustom}, where a transistor sizing problem was addressed by means of optimization techniques. We shall summarise key ideas of the SSTA research in this subsection.

The research at that stage was based on variants of the idea first presented by Clark~\cite{clark1961GreatestFiniteSet} within the block-based approach. The idea is that actual distributions can be approximated by Gaussians by matching the first two moments (mean and variance). Thus, in~\cite{chang2003StatisticalTimingAnalysis,chang2005StatisticalTimingAnalysis} Clark's algorithms were accompanied by handling spatial correlations using principal component analysis (PCA).
To propagate a delay through the timing graph, the linear \emph{canonical model} of a delay was proposed~\cite{visweswariah2004FirstOrderIncrementalBlockBased,chang2005ParameterizedBlockbasedStatistical,visweswariah2006FirstOrderIncrementalBlockBased}. The delay is described as a linear function of parameter variations:
\begin{equation}\label{eq:canonical_form}
  D = a_0 + \sum\limits_{i=1}^n a_i \Delta X_i + a_{n+1} \Delta R,
\end{equation}
where $a_0$ is the mean or nominal value, $\Delta X_i$ represents the variation of $n$ global sources of variation $X_i$ from their nominal values, $a_i$ are the sensitivities to each of the RVs, and $\Delta R$ is the variation of the independent RV, $R$. Then, the mean and variance of a delay were represented using a concept of \emph{tightness probability} (or \emph{binding probability} in~\cite{jess2003StatisticalTimingParametric,jess2006StatisticalTimingParametric}) $T_A=P(A>B)$, which is the probability that the arrival time $A$ is greater than $B$. The linear approximation for $\max(A,B)$ was proposed.

Also, various extensions to this approach were proposed mainly based on adding non-linear terms to~\eqref{eq:canonical_form}. For example, the quadratic term was introduced in \cite{zhan2005CorrelationawareStatisticalTiming,zhang2005CorrelationpreservedNonGaussianStatistical,zhang2006StatisticalTimingAnalysis}. \cite{chang2005ParameterizedBlockbasedStatistical} proposed to use numerically computed tables to describe the non-linear part of the canonical form. \cite{khandelwal2005GeneralFrameworkAccurate} considered gate delays and arrival times using their Taylor-series expansion. The paper~\cite{ramprasath2016SkewNormalCanonicalModel} discusses another modification of the canonical form~\eqref{eq:canonical_form}, based on the addition of the quadratic term and using skew--normal distributions.
The correlations were considered in~\cite{chang2005StatisticalTimingAnalysis}, as we mentioned above, within the PCA method. Later, in~\cite{singh2008ScalableStatisticalStatic}, it was proposed to transform the set of correlated non-Gaussian variables \textit{via}  an independent component analysis (ICA) into a non-correlated set.
The described canonical delay model suffers from a big disadvantage: it requires approximation of the maximum operation, which is a source of errors that we want to mitigate.

More broadly, optimisation problems appear in various aspects of CAD for VLSI~\cite{korte2008CombinatorialProblemsChip,brenner2009AnalyticalMethodsPlacement,held2011CombinatorialOptimizationVLSI}. For example, the gate sizing problem has received attention in the community for more than 30 years~\cite{berkelaar1990GateSizingMOS,sapatnekar1993ExactSolutionTransistor,jacobs2000GateSizingUsing,rakai2015SizingDigitalCircuits}. One should point out the GP formulations of this optimization problem~\cite{boyd2005DigitalCircuitOptimization,boyd2007TutorialGeometricProgramming,joshi2008EfficientMethodLargeScale,naidu2021ConvexProgrammingSolution}, as the latter plays a significant role in the present study, but we shall not discuss these works in detail, as the gate sizing is out of scope in the present paper. One should note that in these papers, as well as in the above mentioned work~\cite{beece2014TransistorSizingCustom}, timing analysis was performed using the canonical model of a delay.

\section{Statistical Static Timing Analysis: Setting up the Problem} \label{sec:ssta_setup}

This Section revisits the problem of calculating the maximum delay in the SSTA histogram approximation, as captured in Algorithm~\ref{alg:SSTA}. In particular, we will present the exact computation of the maximum and the convolution, without claiming novelty of the presented material.

\SetKwComment{Comment}{/* }{ */}
\begin{algorithm}[t!]
\caption{General SSTA algorithm}\label{alg:SSTA}
\KwData{Distributions of delays for $N$ gates and for the input signals}
\KwResult{A distribution of the $\max$ delay of a circuit}

\For{$i \gets 1,..., N$}{
    $M \gets \max$ of arrival times\;
    $C \gets$ convolution of the $i^{\text{th}}$ gate and $M$\;
    propagate $C$ further as input PDF\;
}
$D \gets \max$ of output distributions\; 
\end{algorithm}

Throughout this paper, distributions are represented by the histogram approximation, where we assume that all histograms share edges of all bins. This assumption is made for the sake of clarity and can be removed at the cost of a somewhat more complicated notation.

Let us find the histogram approximation of (i) a distribution of the maximum $\zeta$ of two independent random variables $\eta$, $\xi$, and (ii) a convolution of two histograms. We assume a set of bins $B=\{0,~1,~...,~n-1\}$, and that the histogram samples $\eta_i$ and $\xi_i$ take the value in the interval~$[a, b]$, $\forall{i}$. Given the edges interval $[a, b]$, we partition $\mathbb{R}$ into \textcolor{black}{$n$ intervals of size $|b - a| / n$ with points $e_1$,~$e_2$,~...~$e_{n+1}$ such that $e_i$ and $e_{i+1}$ are the start and end of the $i^{\text{th}}$ interval correspondingly.} The midpoints of the intervals $m_1$,~$m_2$,~...~$m_N$ are also given.

\subsection{Maximum operation}
\label{sec:maximum}

The maximum of two RVs, $\zeta=\max(\eta,\xi)$, which is an RV itself, is defined as 
\begin{equation}\label{eq:max_of_two_rv_def}
    \zeta = 
\begin{cases}
    \eta,& \text{if } \eta \geq \xi,\\
    \xi,& \text{if } \eta < \xi.
\end{cases}
\end{equation}
Using the law of total probability and taking into account the independence of $\xi$ and $\eta$, the probability $P$ of a realization $\zeta=z$ can be written as
\begin{equation}\label{eq:max_of_two_rv_prob}
    P(\zeta = z) = P(\eta = z) \cdot P(\xi \leq z) + P(\xi = z) \cdot P(\eta < z).
\end{equation}
It is easy to see that the discrete random variable formulation and histogram estimations $h_{\zeta}$, $h_{\eta}$, and $h_{\xi}$ for the RVs $\zeta$, $\eta$, and $\xi$ are expressed as follows:
\begin{equation}\label{eq:max_of_two_rv_hist}
    h_{\zeta}[i] = h_{\eta}[i] \cdot \sum_{k = 1}^{i}{h_{\xi}[k]} + h_{\xi}[i] \cdot \sum_{k = 1}^{i-1}{h_{\eta}[k]}.
\end{equation}
Note that the upper bound of the second sum is $i-1$. This is due to a strict inequality, $\eta<\xi$, in \eqref{eq:max_of_two_rv_def}. This is summarized as the Algorithm~\ref{alg:MaxQuad}.


\SetKwComment{Comment}{/* }{ */}
\begin{algorithm}[t!]
\caption{Maximum of two histograms}\label{alg:MaxQuad}
\KwData{Number of bins $n$, two vectors  $x$, $y$ of histogram values}
\KwResult{A new vector $z=\max(x,y)$}
$z \gets \mathbf{0}$\;
\For{$i \gets 1,..., n$ }{
    \For{$k \gets 1,..., i$}{
        $z[i] \gets z[i] + x[i]\cdot y[k]$ \Comment*[r]{first term in \eqref{eq:max_of_two_rv_hist}}
        \If{i != k}{
         $z[i] \gets z[i] + y[i]\cdot x[k]$ \Comment*[r]{second term in \eqref{eq:max_of_two_rv_hist}}
        }
    }

}
\end{algorithm}

\subsection{Convolution operation}
\label{sec:convolution}

The convolution of two discrete-valued functions, $f$ and $g$, is defined as :
\begin{equation}\label{eq:convolution_def}
    (f*g)[z] = \sum_{k=-\infty}^{\infty}{f[k] \cdot g[z-k]}.
\end{equation}
Time complexity of the na\"ive implementation of the convolution is $\mathcal O(N^2)$, which is can be seen in Algorithm~\ref{alg:convolution_exact}. The formula \eqref{eq:convolution_def} implies that the values of the edges must be changed. The value of the first edge has to be added to all other edges. This can be done in many ways and is discussed in the following.

One can add the first value to all edges and unite them during the SSTA algorithm when the edges of the second histogram differ. The receipt is simple: find a new array of edges $\mathbf{e}$ and modify the PDFs of histogram approximations $f_{\alpha}$, $f_{\beta}$ with the new changed edges. The array of edges of $f_{\alpha}$ is given as $\mathbf{e}^{\alpha}$, similarly $\mathbf{e}^{\beta}$ denotes the array of edges of $f_{\beta}$.


Similar to the \texttt{rv\textunderscore histogram()} method of the \texttt{scipy.stats} library, one can find a distribution function (PDF) that fits the given histogram as follows. Let us say that $F$ is the fitted cumulative distribution function (CDF) of the $f_{\alpha}$. Then PDF of the realization $z \in n_i$ of the new histogram $f_{\alpha'}$ with the desired edges $\mathbf{e}$ is
\begin{equation}\label{convolution_integral}
    f_{\alpha'}(z) = F(e_{i+1}) - F(e_{i}),
\end{equation}
or, recalling a definition of CDF,
\begin{equation}\label{convolution_union_integral}
    f_{\alpha'}(z) = \int_{e_i}^{e_{i+1}} f_{\alpha}(x) \,dx.
\end{equation}
The same would be done for $f_{\beta'}$. The solution \eqref{convolution_union_integral} gives more precise results and bypasses the problem of fitting functions, which is relatively computationally demanding, all at the cost of a slightly longer code. The exact integration can be performed in $O(N)$ for the whole histogram.

\begin{figure}[t]
\centering
\includegraphics[width=\textwidth]{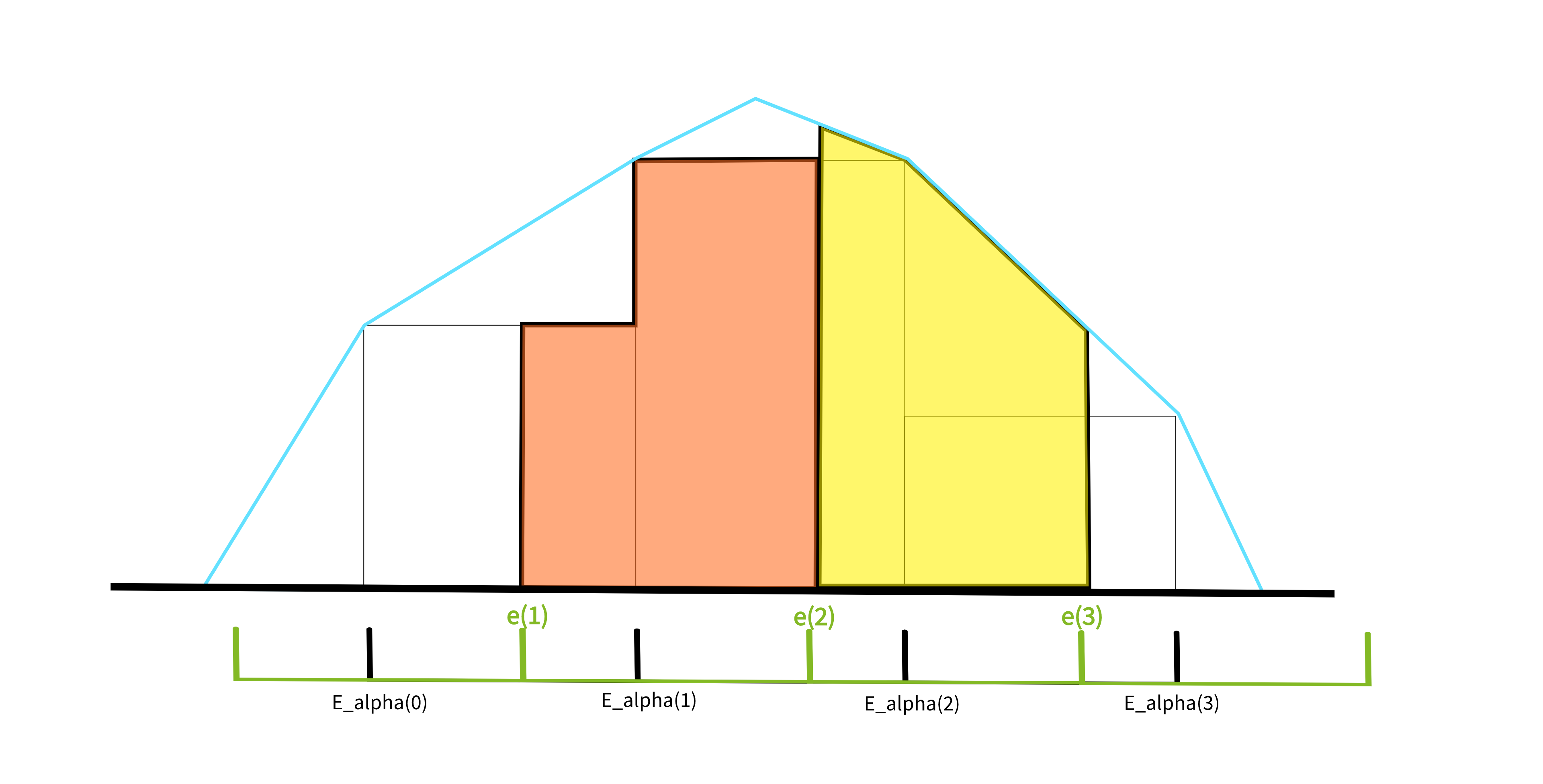}
\caption{Sketch demonstrating integration over the fitted function \eqref{convolution_integral} (yellow surface) and the exact integration \eqref{convolution_union_integral} (red surface). Blue function represents the fitted PDF.}
\label{sketch}
\end{figure}

When one looks at the SSTA Algorithm~\ref{alg:SSTA}, a problem with such a union of edges is evident. After convolutions of the input gates, every time a maximum and then convolution are computed, the edges will differ and have to be united. Taking into account the two inputs for each gate, the function \eqref{convolution_integral} or \eqref{convolution_union_integral} is called twice per a gate. Moreover, more problems are to come when trying to solve a convolution optimization problem.

A different and more straightforward solution is presented in the next subsection.

\subsection{Convolution with shifting a histogram}\label{hist_shift}

A second solution to the problem of adding a value to the edges is a simple shifting. We can shift the whole histogram to the left or right by the number of bins determined by a value that is to be added to all edges divided by the length of the bins and floored. Shifting the value of a bin by such a number simulates the addition of the first value to all edges. We assume $n$ bins, set of bins $B = \{0, 1, ..., n - 1\}$ the identical edges of the histograms is given as $\mathbf{e} \in \mathbb{R}^{(n+1)\times 1}$, two RVs  $\alpha, \beta$; their convolution $\zeta$, its shifted version $\zeta'$; and their histogram approximations $h_{\alpha}$, $h_{\beta}$, $h_{\zeta}$, $h_{\zeta'}$. The shift $s$ can be computed as 

\begin{equation}
s = \left\lfloor \frac{\mid e_{0}\mid}{e_{1} - e_{0}} \right\rfloor,
\end{equation}
where $\lfloor \cdot \rfloor$ denotes the flooring operation.
In the case of $e_{0} > 0$, the new changed histogram $h_{\zeta'}$ at point $x \in B$ will look like
\begin{equation}
h_{\zeta'}[x + s] = h_{\zeta}[x]
\end{equation}
In the case of $e_{0} < 0$, the shift is similar: 
\begin{equation}
h_{\zeta'}[x] = h_{\zeta}[x + s]
\end{equation} 

When shifting to the right, there are $s$ unoccupied positions on the left. These are nullified. Similarly, done when shifting to the left. Having the starting interval set correctly, this does not have any effect on precision as the starting and ending bins should always be zero. If the interval is small, the accuracy increases, since the bins can encode more information in a smaller interval. However, if it is too small, then we cut some information by this shift. Furthermore, the more bins we add, the more precise this shift will be.

\SetKwComment{tcp}{// }{}%
\begin{algorithm}[t!]
\caption{Convolution}\label{alg:convolution_exact}
\KwData{Number of bins $n$, two vectors  $x$, $y$ of histogram values, edges array $e$}
\KwResult{Convolution of two histograms in $c$}\label{alg:convolution}
$c \gets \mathbf{0}$\;
       \tcp{perform convolution}
\For{$z \gets 1,..., n$}{
    \For{$k \gets 1,..., z$}
    {$c[z] \gets c[z] + x[k]\cdot y[z-k]$\;}
}
    \tcp{shift histogram}
$s \gets \text{floor}(\text{abs}(e[0])/(e[1]-e[0]))$

\If{$e[0] > 0$}{
$c[s:] \gets c[:-s];$\\
$c[:s] \gets 0;$
}
\If{$e[0] < 0$}{
$c[:-s] \gets c[s:];$\\
$c[-s:] \gets 0;$
}
\end{algorithm}

The shifting method gives exactly precise solutions as the one with union. It is \textcolor{black}{straightforward} to implement, can be done in linear time, and can be used nearly without any change in the optimization problem. Therefore, this method is used better than the union method. 

{\color{black}
In summary, we have reviewed two key (atomic) operations of SSTA algorithms, maximum and convolution, and the way how these operations can be performed for histograms. In the following sections, we will give a formulation of the SSTA as (i) a Binary Integer Programming problem and (ii) a Geometric Programming problem. These formulations given for histograms constitute the original contribution of this work.
}
\section{SSTA via Binary--Integer Programming}
\label{sec:ssta_bip}

Here, we formulate the SSTA algorithm as a Binary--Integer Programming (BIP) problem. For this purpose, we will introduce a unary encoding of {\color{black} counts in binary variables\footnote{\color{black}This may seem confusing, at first, but we do utilize \emph{binary} decision variables to store \emph{unary}-encoded counts. The \emph{unary} encoding allows for an easier formulation of the corresponding constraints than in the case of \emph{binary} encoding in \emph{binary} variables, which is also possible \cite{vielma2011modeling}.}}. Then, we will show how to perform the atomic operations ($\max$ and convolution) on histograms \textit{via} BIP. Finally, we will compare the approach against Monte Carlo simulations and discuss its scalability.

\subsection{Statement of the Problem}\label{sec:ssta_bip:problem_statement}

To write the SSTA Problem of a Circuit, $D_{\max} = \max(D_1,\ldots,D_N)$, as a BIP problem \eqref{eq:MIP_general}--\eqref{eq:BIP_constraints_general}, one should (i) propose a risk measure for the objective function $\mathbf{c}^T \mathbf{x}$, and (ii) formulate corresponding constraints. Let us discuss the constraints first.

From the Algorithms~\ref{alg:MaxQuad} and \ref{alg:convolution}, one can see that a multiplication of two non-negative real numbers occurs in both of them. By utilizing multiplication naively, we obtain a bi-linear function, which is not convex jointly in both arguments. This cannot be used as a constraint or as an objective function. One solution could be to use McCormick envelopes. This requires setting the lower and upper bounds of the factors. The only way to compute these bounds is by the exact computation of the problem using the methods shown in Section~\ref{sec:ssta_setup}. Another option is to formulate the problem in unary notation\textcolor{black}{, which is the subject of the present Section}.

{\color{black}

Below, we will discuss how the atomic operations of SSTA can be performed on histograms in unary encoding. It will be shown that this leads to the corresponding mixed--integer linear programs.


\subsection{Atomic operations on histograms in {\color{black} unary} representation}
\label{sec:ssta_bip:atomic}

As we have discussed above, multiplication is the key operation for both $\max$ and convolution. This Section discusses the multiplication in {\color{black} unary} representation.

\subsubsection{Multiplication}

Let us first make a note on vectorization of the multiplication operation.
Consider two natural numbers, $\alpha$ and $\beta$. By definition, the multiplication of two numbers is equal to the repeated addition:
\begin{subequations}
\begin{equation}
    \alpha \cdot \beta = \underbrace{\beta + \beta + \ldots + \beta}_{\text{$\alpha$ times}}.
\end{equation}
It is easy to see that in unary notation it can be represented by matrix--vector multiplication. Indeed, having written $\alpha$ and $\beta$ as column vectors, $\mathbf{a}$ and $\mathbf{b}$, of sizes $\alpha\times1$ and $\beta\times1$, correspondingly, we obtain
\begin{equation}\label{matrixMult}
\alpha \cdot \beta =
\underset{1\times \alpha}{\mathbf{1}^T}
\underset{\alpha \times \beta}{
    \begin{bmatrix}
    \mathbf{a} & \mathbf{a} & \ldots & \mathbf{a}\\
    \end{bmatrix}
}
\underset{\beta \times 1}{\mathbf{b}}
\end{equation}
\end{subequations}
This representation allows performing multiplication of numbers written in unary encoding efficiently. Now, we shall discuss multiplication of \emph{binary variables} by means of integer programming.

Consider two binary variables: $x,y\in\{0,1\}$. Multiplication of the variables via BIP requires introduction of an auxiliary variable, $s\in\{0,1\}$, and setting corresponding constraints. The constraints establish obvious relationship between this variable, $s$, and the multiplicands, $x$ and $y$:
\begin{equation}\label{binary_constr}
\left.
\begin{aligned}
s &\leq x,\\ 
s &\leq y,\\ 
s &\geq x + y - 1.\quad
\end{aligned}
\right\}
\end{equation}


The process of formulating the constraints \eqref{binary_constr} for the convolution of two histograms in binary form is summarized as Algorithm~\ref{alg:convolution_bip}. Note that the two outer loops, over $z$ and $k$, correspond to those from Algorithm~\ref{alg:convolution_exact}, and the inner loops, over $i$ and $j$, are due to {\color{black} unary encoding of the counts}. However, these inner loops can be removed if the computation performed via efficient vector operations as shown in \eqref{matrixMult}.

The BIP constraints for the $\max$ operation can be introduced in a similar manner and, hence, not discussed here. This procedure will repeat two loops from Algorithm~\ref{alg:MaxQuad} and will have two internal loops as in Algorithm~\ref{alg:convolution_bip}. Therefore, further in this Section we will discuss only convolution operation, but we want to note that the same argumentation can be carried out for the maximum of two histograms.

\SetKwComment{Comment}{/* }{ */}
\begin{algorithm}[t!]
\DontPrintSemicolon
\caption{Generation of BIP constraints for convolution operation 
}\label{alg:convolution_bip}
\KwData{Number of bins $n$, number of binary variables $m$, two histograms in {\color{black} unary encoding}, $\mathbf H_X$ and $\mathbf H_Y$, of size $n\times m$}
\KwResult{1-D array $\mathbf u$ of size $n$ with auxiliary variables $s$ and a list with the corresponding constraints}
$\mathbf u \gets \mathbf{0}$\;
\tcp{\textcolor{black}{for all bins of a histogram $\mathbf H_X$}}
\For{$z \gets 1,\ldots, n$}{
    \tcp{\textcolor{black}{for all bins of a histogram $\mathbf H_Y$}}
    \For{$k \gets 1,\ldots, z$}
    {
        \tcp{\textcolor{black}{for all unary digits of the bin of the histogram $\mathbf H_X$}}
        \For{$i \gets 1,\ldots, m$}
        {
            \tcp{\textcolor{black}{for all unary digits of the bin of the histogram $\mathbf H_Y$}}
            \For{$j \gets 1,\dots, m$}
            {
                initialize a variable, $s$\;
                $\mathbf u[z] \gets \mathbf u[z] + s$\;
                $x \gets \mathbf H_X[k,i]$\; 
                $y \gets \mathbf H_Y[z-k,j]$\;
                \textcolor{black}{add unary constraints \eqref{binary_constr} linking $s$, $x$, and $y$ for the specific pair $(i,j) $ of unary digits in the histograms $\mathbf H_X, \mathbf H_Y$ }\;
            
            }
         }
    }
}
\end{algorithm}

\subsubsection{Forming a histogram in the {\color{black} unary encoding}}

Consider two histograms in {\color{black} unary encoding}, $\mathbf H_X$ and $\mathbf H_Y$, of size $n\times m$, where $n$ is the number of bins and $m$ is the number of binary variables in a bin. For these histograms, the Algorithm~\ref{alg:convolution_bip} returns an array $\mathbf u$ with auxiliary variables $s$ given as variables. The values of these variables shall be obtained during the solution of the BIP problem.

One can see that the length of the array $\mathbf u$ equals the number of bins, $n$. Each element of the array $\mathbf u$ contains a sum of the auxiliary variables (see Algorithm~\ref{alg:convolution_bip}, lines~$6-7$). It is easy to see the meaning of the each element in $\mathbf u$ (and, hence, of the each sum of $s$ variables). They correspond to the values of bins in a histogram $\mathbf H_{XY}$, which is the  convolution of $\mathbf H_X$ and $\mathbf H_Y$.

In order to map the elements of $\mathbf u$ onto the bin values of the histogram $\mathbf H_{XY}$, we (i) generate the matrix $\mathbf H_{XY} \in \{0, 1\}^{n\times m}$ with \emph{new} variables, and (ii) supplement these variables with the following constraints
\begin{equation}\label{eq:constr_bip_hist_norm}
    \left.
    \begin{aligned}
    \mathbf 1^T \mathbf H_{XY}[1,:]  &\leq {\color{black}\mathbf u[1]} \frac1d + 0.5 \\
    \vdots \\
    \mathbf 1^T \mathbf H_{XY}[n,:]  &\leq {\color{black}\mathbf u[n]} \frac1d + 0.5
    \end{aligned}
    \quad
    \right\}
\end{equation}
Here $\mathbf 1^T \mathbf H_{XY}[i,:]$ denotes summation over the $i^{\text{th}}$ row of the histogram $\mathbf H_{XY}$, \textit{i.e.} $\mathbf 1^T \mathbf H_{XY}[i,:] = \sum_{k=1}^m \mathbf H_{XY}[i,k]$, and {\color{black}$\mathbf u[i]$} is the $i^{\text{th}}$ element of the array $\mathbf u$. Recall that the histogram $\mathbf H_{XY}$ has $n$ rows and $m$ columns, where rows correspond to the bins and columns give the bin counts in {\color{black} unary encoding}. The parameter $d$ is a normalization factor that reads
\begin{equation}\label{eq:bip_normalization_factor}
    d = \frac{
    \max\{
    {\color{black}\mathbf u[1], \ldots, \mathbf u[n]}\}
    }
    {m}.
\end{equation}
Such a choice of the normalization parameter is to ensure that the r.h.s. of the inequalities~\eqref{eq:constr_bip_hist_norm} does not exceed the number of binary variables, $m$. The number $0.5$ is to round the r.h.s. to a positive real number.

In this form the value of the normalization factor $d$ can be obtained self-consistently during the optimization procedure. \textcolor{black}{However, this would give a non-convex optimization problem due to the division in \eqref{eq:bip_normalization_factor}, which we wish to avoid. Therefore, we fix the normalization factor in our implementation to be a constant of our choice.} We can check post-hoc, whether the solution has reached the bound on the number of unary digits in any of the bins of the histogram. If this is the case, the optimality has been affected and we can increase the bound and re-run the procedure. If this is not the case, the normalization does not affect the optimality.

Let us now discuss how the problem can be simplified by the tightening.

\subsubsection{Problem tightening}
The convergence of the BIP solver can be increased by introducing constraints that tighten the relaxation. For example, one can separate zeroes from ones in the histogram matrix $\mathbf H$ with symmetry--breaking constraints:
\begin{equation}\label{eq:constr_symmetry-break}
    H_{i, 1} \geq H_{i, 2} \geq H_{i, 3} \geq \ldots \geq H_{i, m-1} \geq H_{i, m}.
\end{equation}
These constraints do not have any effect on the correctness of the solution but decrease size of the branch and bound tree of the solver.

Moreover, as we discussed in Section~\ref{sec:background:definitions} (see expressions \eqref{eq:gate_operation} and \eqref{eq:gate_operation_rv}), maximum and convolution describe the operation of a basic logic gate. Thus, it is natural to evaluate these operations simultaneously for each gate. Similar to \eqref{binary_constr}, we can write for three binary variables, $x,y,z\in\{0,1\}$:
\begin{equation}\label{eq:const_3-term_multiplication}
\left.
    \begin{aligned}
    s &\leq x,\\ 
    s &\leq y,\\
    s &\leq z,\\ 
    s &\geq x + y + z - 2.
    \end{aligned}
    \right\}
\end{equation}
These constraints describe the $\max$ of two binary variables, $x$ and $y$, and further convolution of the result with the third binary variable, $z$. 
Having discussed the atomic SSTA operations in unary encoding, we can proceed to the formulation of the SSTA Problem as a BIP problem,
{\color{black} where we perform the same operation bit-wise on the unary-encoded counts.}

\subsection{Implementation and validation}\label{sec:ssta_bip:implementation}

Calculation of the circuit delay requires traversing the timing graph $G(E,V)$. Since a histogram $\mathbf H_{\text{sink}}$ corresponding to the sink contains all binary variables obtained in the previous steps, it is natural to write the objective function of the BIP problem \eqref{eq:MIP_general}--\eqref{eq:BIP_constraints_general} as the sum of all the variables due to atomic operations in the graph, $\mathbf{1}^T \mathbf H_{\text{sink}} \mathbf{1}$, subject to constraints discussed above. Doing so, we obtain the SSTA problem of a Circuit \textit{via} BIP as follows:
\begin{alignat*}{3}\label{eq:ssta_bip_final}\numberthis
 & \text{minimize} & {\color{black} \textrm{risk}(\mathbf H_{\text{sink}})} &\\ 
 & \text{subject to} 
                \quad & \mathbf H_g[i,1] \geq \ldots &\geq \mathbf H_g[i,m]  \quad & \forall g \in G: i&=1,\ldots,n,\\
                & & \mathbf 1^T \mathbf H_g[i,:]  &\leq {\color{black}\mathbf u[i]} \frac1d + 0.5, & & \\
                & & s_g &\leq \{x,y,z\},  & &\\
                & & s_g &\geq x + y + z - 2,  & &\\
                & & \mathbf{S}_g &\leq \mathbf{G}_g   &  & 
\end{alignat*}
{\color{black} where $\textrm{risk}$ is a MILP representation of a risk measure. Notice that in many practical scenarios, one may consider the conditional value at risk (CVaR) \cite{rockafellar2000optimization} of the histogram--approximated random variable 
$\mathbf H_{\text{sink}}$ as the risk measure $\mathrm{risk}(\mathbf H_\text{sink})$, in an effort to ``shift the probability mass left'', loosely speaking. In general, this would be implemented by summing up the counts in some number of ``rightmost'' bins of the histogram approximation. For testing purposes, we have utilized the expression $\mathbf{1}^T \mathbf H_{\text{sink}} \mathbf{1}$, which sums counts across
the histogram, and implements CVaR at 0\% confidence level}.

\begin{figure}[t]
\centering
\includegraphics[width=0.95\textwidth]{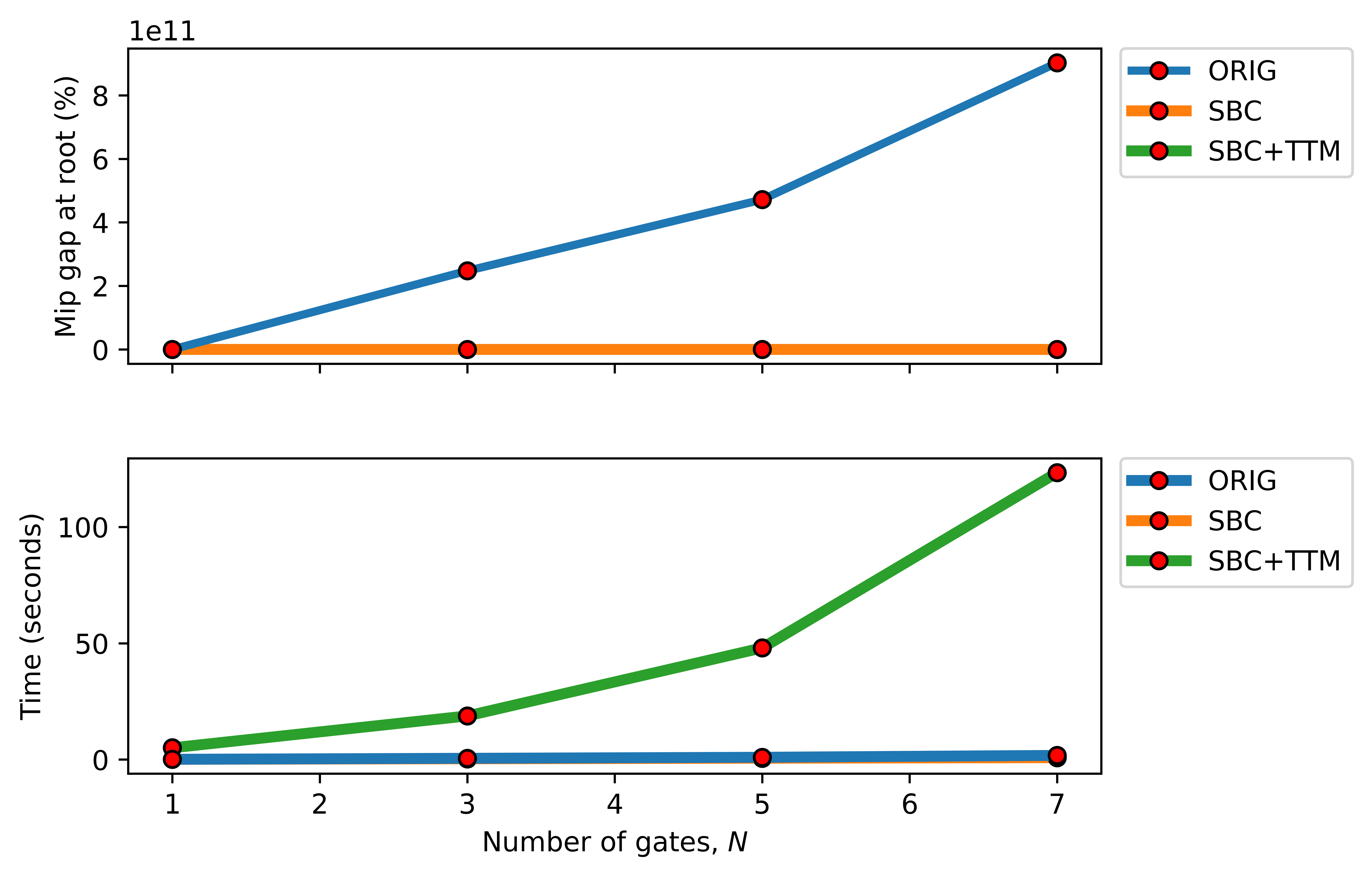}
\caption{Comparison of three BIP formulations for the SSTA. The methods are tested on a ``ladder'' of maxima with $n=12$ bins, {\color{black} each bins' count bounded from above by} $m=12$. The blue line indicates a method with only \eqref{binary_constr} constraints  (ORIG), the orange line indicates a method with the symmetry--breaking constraints (SBC) \eqref{eq:constr_symmetry-break} and the green line indicates a method with the symmetry--breaking constraints and a $3-$term multiplication model (SBC+TTM)\eqref{eq:const_3-term_multiplication}. In the first figure, the orange line overlaps the green line. In the second figure, the blue line overlaps the orange line.}
\label{relax_comparison}
\end{figure}

This BIP formulation uses the symmetry--breaking constraints~\eqref{eq:constr_symmetry-break} and the $3-$term multiplication model~\eqref{eq:const_3-term_multiplication}. Note that these constraints \eqref{eq:const_3-term_multiplication} imply taking $\max$ and convolution operations for a gate $g$, whereas constraints \eqref{binary_constr} correspond only to multiplication and, thus, should be added after each atomic operation, both $\max$ and convolution. Matrices $\mathbf{G}_g$ in the BIP are unary representations of a delay of each gate, bounded from above elementwise by $\mathbf{S}_g$, again represented in unary.

\begin{figure}[t]
\centering
\includegraphics[width=\textwidth]{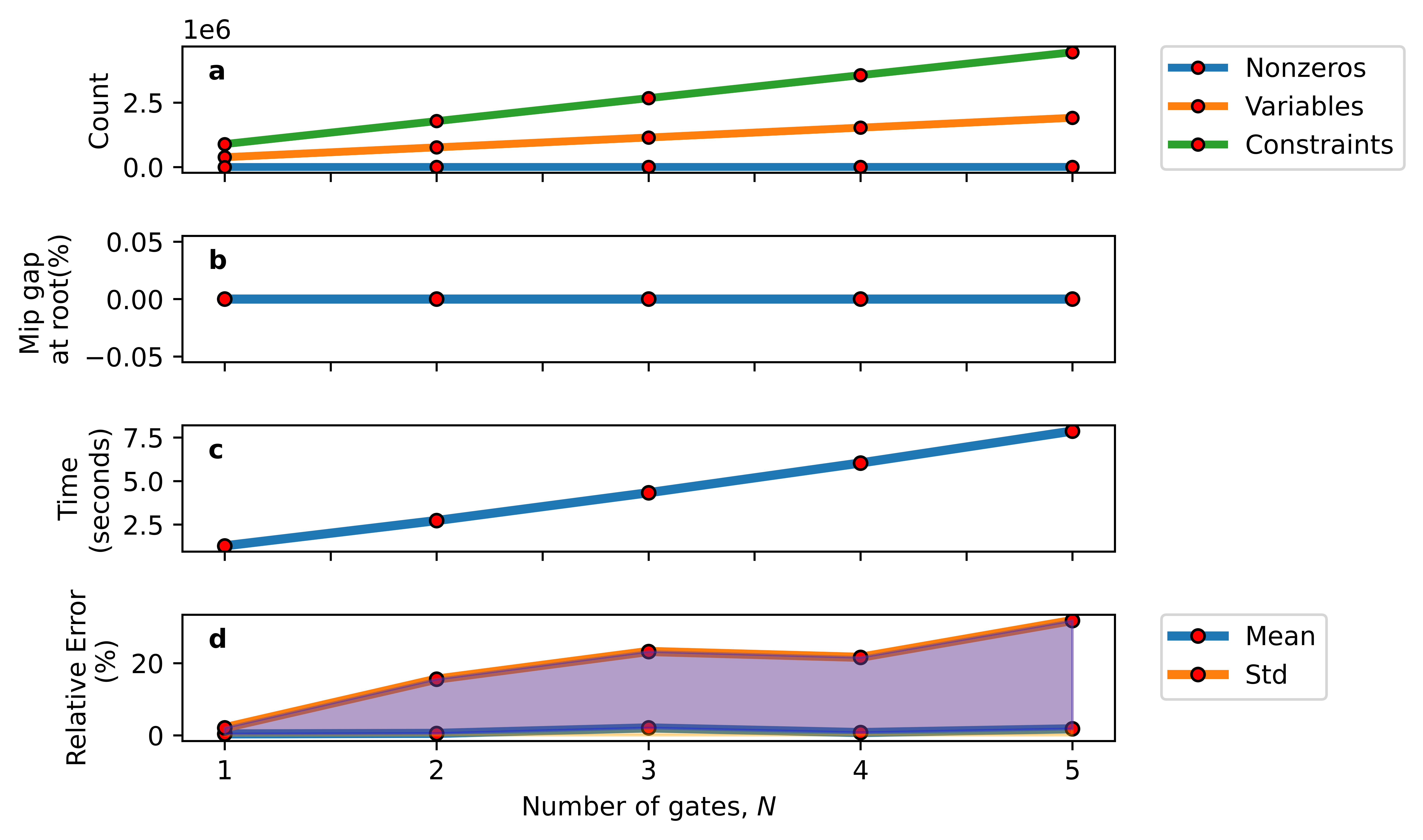}
\caption[Scalability of the model with the first relaxation]{Scalability of the model with the first relaxation (SBC) \eqref{eq:constr_symmetry-break} tested on a ``ladder'' of maxima with $n=25$ bins, $m=20$ unary variables, and no time limit. The subplots show:  \textbf{a},~ the growth in the number of non-zeros (blue line), variables (orange line), and constraints (green line). \textbf{b},~MIP gap at a root node in percentage, MIP gap tolerance is set to 1\%; \textbf{c},~time in seconds; \textbf{d},~relative error of the standard deviation (orange line) and mean (blue line) compared to Monte Carlo.}
\label{scaling}
\end{figure}

To validate our BIP formulation of the SSTA, we have implemented \eqref{eq:ssta_bip_final} in 
Mosek~10.0 matrix-oriented API. As a test bench we have used a toy circuit used in \cite[Fig.~7]{mishagli2020RadialBasisFunctions} that gives a ``ladder'' sequence of logic gates. In terms of atomic operations, this sequence for the $N$th gate reads
\begin{equation}\label{eq:ladder_sequene}
     \max\{\underbrace{\ldots\max[\max(\xi_1,\xi_2) + \xi_0, \xi_3] + \xi_0\ldots}_{N-1\text{ times}},  \xi_{N+1}\} + \xi_0.
 \end{equation}
Here RVs $\xi_i$ are drawn from the normal distribution, $\xi_i\sim\mathcal N(\mu_i,\sigma_i)$. Delay due to gates operation time (gate delay) is given by $\xi_0$, and $\xi_i$ ($i=1,\ldots,N+1$) are inputs' delays. Gate delays were assumed distributed according to the standard normal distribution, $\xi_0\sim\mathcal N(1,0)$; the mean values and standard deviations for inputs were drawn from the uniform distribution, following Ref.~\cite{mishagli2020RadialBasisFunctions}. 

This sequence was simulated (i) using Monte Carlo (MC) and (ii) by solving BIP \eqref{eq:ssta_bip_final}. For each gate $g \in G(E,V)$, a new matrix $\mathbf H_g$ containing binary variables was created as described above, until the sink was reached. Then, the final BIP problem was passed to a Mosek solver. Numerical experiments were ran on a machine equipped with Intel(R) Core(TM) i9-9880H (8 cores at 2.3~GHz and total of 16 threads) with 16~GB RAM. 

The results are summarized in Figs.~\ref{relax_comparison} and \ref{scaling}. Correctness of the BIP formulation can be seen from Fig.~\ref{scaling} where the mean and standard deviation of the sink delay distribution is compared against the MC. Although the number of constraints and variables scale linearly with the number of gates, these numbers are huge, thus, this approach does not scale. For example, the numbers of variables and constraints quickly reached the order of $\sim10^6$ for only 5 gates (Fig.~\ref{scaling}). In the next Section, we consider the GP formulation of SSTA, which scales much better.
}

\section{SSTA via Geometric Programming}
\label{sec:ssta_gp}

Next, we present \textcolor{black}{a practical} formulation \textcolor{black}{of the SSTA problem}.
First, we will present \textcolor{black}{the formulation by means of Geometric Programming}, which is a restriction of the exact formulation, but where the error can be made arbitrarily small. Then, we will show its reformulation \textcolor{black}{and} scalability.

{\color{black}
\subsection{Statement of the Problem}}
\label{sec:ssta_gp:problem_statement}

\textcolor{black}{Naturally}, we can treat the probability of each bin as a positive number in the range $[\epsilon, 1]$, where $\epsilon$ is a very small number. Multiplication of two bins leads to a monomial function of two variables and a neutral coefficient. Either the convolution or the maximum is then the sum of the multiplications, thus a posynomial. \textcolor{black}{Therefore, general Algorithms \ref{alg:MaxQuad} and \ref{alg:convolution} can be utilized within the Geometric Programming (GP) framework in a straightforward manner and no additional constraints are needed for the atomic operations.}

In the following, we again consider \textcolor{black}{a timing graph $G(E,V)$ consisting of $N$ gates; the number of bins used for the histogram approximation of distributions is $n$}. For each input gate $g$, we have a vector $\mathbf{e}_g \in \mathbb{R}_{++}^{n\times 1}$ created from generated numbers with Gaussian probability and a vector $\mathbf{z}_g \in \mathbb{R}_{++}^{n\times 1}$ of non-negative variables representing the bin probabilities. 
Similarly to Section~\ref{sec:ssta_bip}, the geometric program \textcolor{black}{starts with}
\begin{equation}
\begin{aligned}\label{eq:GP_SSTA}
 & \text{minimize} & {\color{black} \textrm{risk}(\mathbf z_{\text{sink}})} \\
 & \text{subject to} \quad & \mathbf{e}_g \leq \mathbf{z}_g & \leq \mathbf{1}, \quad g = 1,\ldots,N,
\end{aligned}
\end{equation}
{\color{black}
where ${\color{black} \textrm{risk}(\mathbf z_{\text{sink}})}$ is a posynomially--representable risk measure of the delay
$\mathbf z_{\text{sink}}$ at the sink of the graph,
and the bounds of the variables are standard GP-compatible inequalities. 
Subsequently, the construction of the geometric program follows the Algorithm~\ref{alg:SSTA} (``General SSTA algorithm'').
In particular, for each maximum $z_{\zeta}$ of two histograms $z_{\eta}$ and $z_{\xi}$, we constrain $z$ as in Algorithm~\ref{alg:MaxQuad}:
\begin{align}
z_{\zeta}[i] = z_{\eta}[i] \cdot \sum_{k = 1}^{i}{z_{\xi}[k]} + z_{\xi}[i] \cdot \sum_{k = 1}^{i-1}{z_{\eta}[k]} \quad \forall i = 1 \ldots n.
\end{align}
See Section~\ref{sec:maximum} for a discussion. 
For each convolution $z_{\zeta}$ of two histograms $z_{\eta}$ and $z_{\xi}$, we constrain elements of $\mathbf{z}$ as in Algorithm~\ref{alg:convolution_exact}:
\begin{align}
z_{\zeta} = \sum_{i=1}^n \sum_{k=1}^{i-1} z_{\eta}[k] \cdot z_{\xi}[i-k]
\end{align}
up to the shifting. See Section~\ref{sec:convolution} for further details.
}

{\color{black}
\subsection{Reformulation and relaxation}\label{Gp_relax}
}
The posynomial formulation in \eqref{eq:GP_SSTA} does not seem to be problematic in any way. 
Still, to improve its scalability, we may wish to consider a reformulation. 
In the following, we will concentrate only on convolution; the procedure is the same for the maximum. 


After the first convolution in the last bin, we have a posynomial with $1\cdot n$ terms (monomials with two variables). After the second convolution, we will have in the last bin a posynomial with $n\cdot n$ monomials each with three variables, and for the $N^{\text{th}}$ gate we will have $n^{N-1}$ monomials in the last bin, each with $N$ variables. This clearly leads to an exponential growth in monomials after each convolution and maximum for a constant number of bins. 

For each monomial in the posynomial, we need to introduce an exponential cone, two continuous variables, and two constraints. Thus, for a constant number of bins, the number of variables, cones, and constraints grows exponentially with the number of gates. For a constant number of gates $N+1$, the number of variables, cones, and constraints grows as the sum of all bins $\sum_{i=1}^n i^N$ with the number of bins, 
which, in turn, can be expressed as a polynomial of degree $N+1$ by Faulhaber's formula~\cite{knuth1993JohannFaulhaberSums}.

\begin{figure}[t]
\centering
\includegraphics[width=\textwidth]{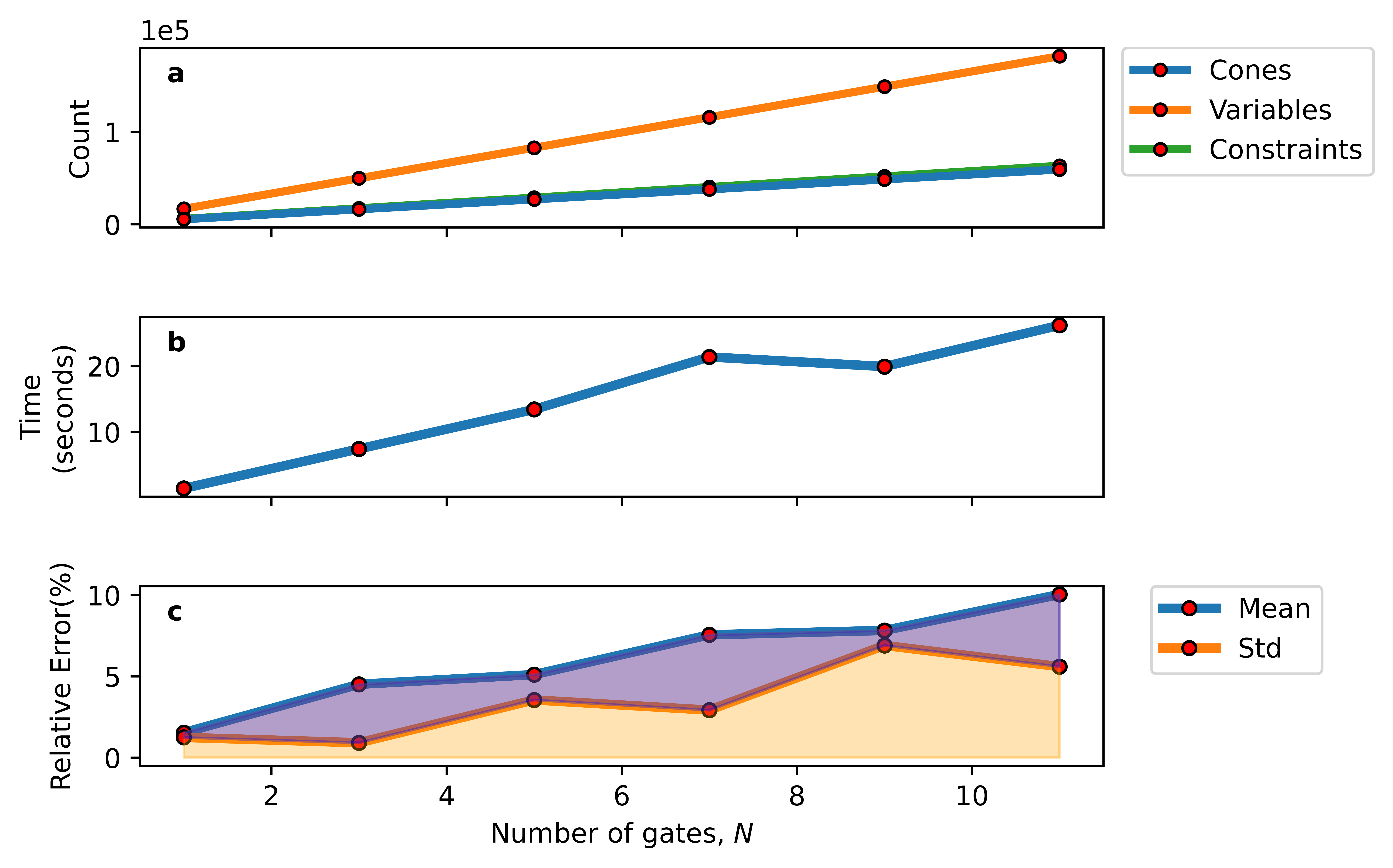}
\caption{Scalability of the GP  \eqref{eq:GP_SSTA} on a ``ladder'' of maxima with $n=60$ bins and a varying number of gates, $N$. The subplots show: \textbf{a},~the growth in the number of cones (blue line), variables (orange line), and constraints (green line overlapped the blue line) as the number of gates, $N$, increases; \textbf{b},~time in seconds; \textbf{c},~relative error of the standard deviation (orange line) and mean (blue line) compared to Monte Carlo as the number of gates, $N$, increases.}
\label{scaling_GP}
\end{figure}

\begin{figure}[t]
\centering
\includegraphics[width=\textwidth]{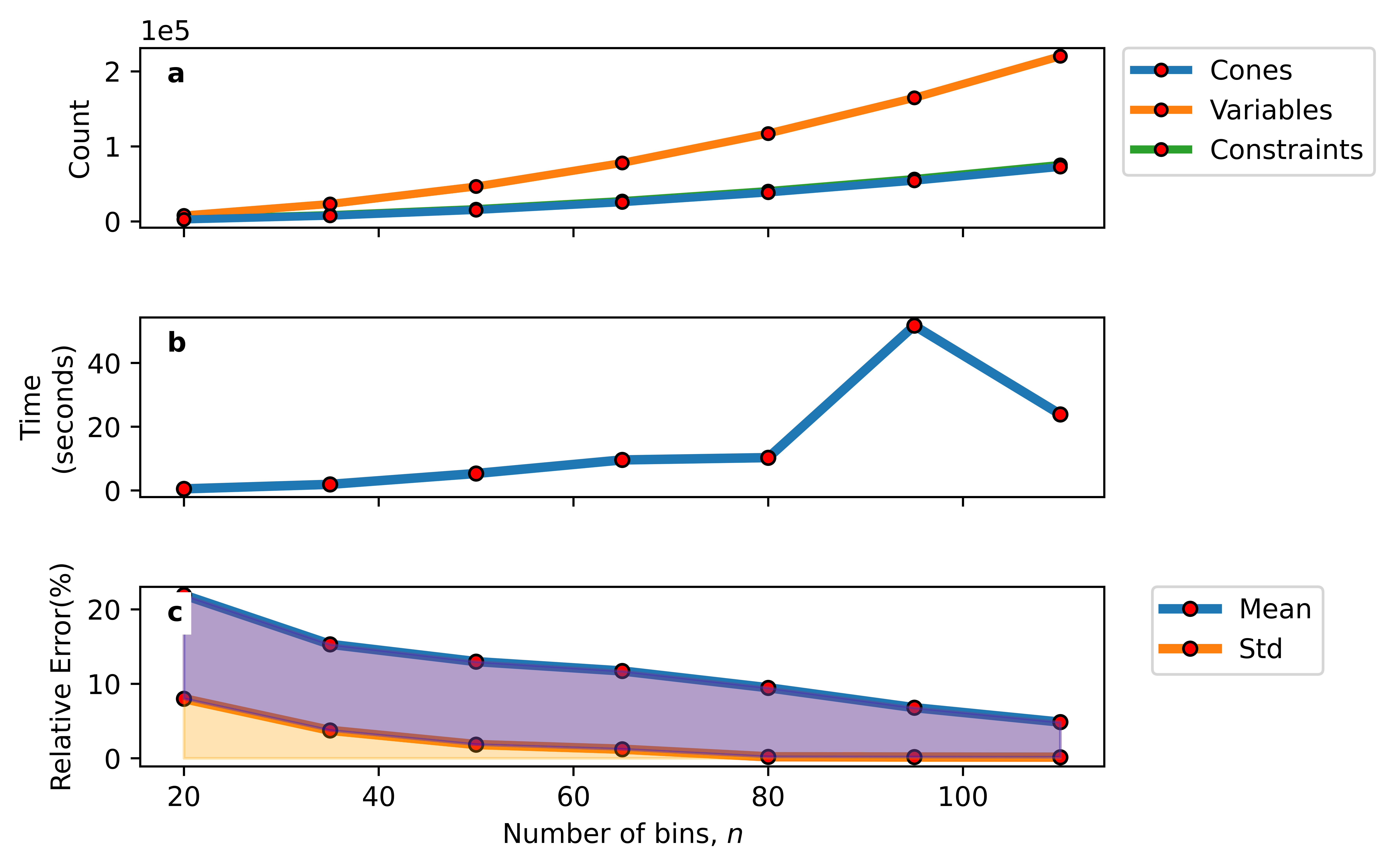}
\caption{Scalability of the reformulated GP \eqref{eq:GP_SSTA} tested on a ``ladder'' of maxima with fixed $N=4$ gates and varying numbers of bins, $n$, per gate. The subplots show: \textbf{a},~the growth in the number of cones (blue line), variables (orange line), and constraints (green line overlapping the blue line) as the number of bins, $n$, increases; \textbf{b},~time in seconds as the number of bins increases; \textbf{c},~relative error of the standard deviation (orange line) and mean (blue line) compared to Monte Carlo, as the number of bins, $n$, increases.}
\label{scaling_GP_bins}
\end{figure}

We can reduce this with the following simple trick: by introducing $n$ new positive variables (monomials) and setting appropriate constraints. 
{\color{black}At the beginning of the traverse, we initiate two empty lists of vectors to store the successors and predecessors. After the convolution at the $i^{\text{th}}$ gate, the resulting vector of posynomials $\mathbf z_i$ is appended to the list of predecessors.  A new vector of one-variable monomials is created and written in the successors' list, and constrained such that it is no less than the predecessor. This new vector represents a histogram and is propagated further. The last vector $\mathbf z_{\text{sink}}$ that appears in the successors' list corresponds to the histogram of the delay at the sink node.
}

Subsequently, the formulation continues with the inequalities based on the equalities in Algorithm~\ref{alg:MaxQuad} and Algorithm~\ref{alg:convolution_exact},
while adding the auxiliary variables.\footnote{For the implementation, see \texttt{maximum\_GP\_OPT} and \texttt{convolution\_GP\_OPT} in \url{https://github.com/bosakad/SSTA-via-GP/blob/experimental/src/timing/cvxpyVariable.py}.}
In particular, for each convolution, we introduce $(n/2)(1 + n)$ new exponential cones, thus $(2n/2)(1 + n) + n$ auxiliary variables.
Such a reformulation gives the exact same solution as the original GP.
We just decreased the exponential growth of variables, cones and constraints to a linear one with the number of gates, and a high-degree polynomial growth with the number of bins to always quadratic. The values are slightly different for the maximum, but the asymptotics before and after the reformulation are identical.

 \color{black} Notice that \eqref{eq:GP_SSTA} can be transformed into a standard GP-compatible inequality.
Thus, the reformulation of \eqref{eq:GP_SSTA} remains a generalized GP program. \color{black}

{\color{black}
\subsection{Implementation and validation}
}
\label{sec:ssta_gp:implementation}

We have prototyped the formulations in CVXPY \cite{agrawal2018RewritingSystemConvex,diamond2016CVXPYPythonembeddedModeling}. For benchmarking purposes, we have passed the instances to MOSEK~10.0, which ran on a laptop equipped with Intel(R) Core(TM) i9-9880H (8 cores at 2.3GHz and total of 16 threads) with 16~GB RAM. \textcolor{black}{The same toy circuit was used as in Section~\ref{sec:ssta_bip:implementation}.}

The scalability of the reformulated GP model is demonstrated in Fig.~\ref{scaling_GP} and Fig.~\ref{scaling_GP_bins}. The results on a ladder of maxima parameterized by the depth of the ladder are shown in Fig.~\ref{scaling_GP} Notice that the numbers of cones (blue line) and variables (orange line) are linear in the depth of the ladder. 
At the same time, the relative error increases.

Fig. \ref{scaling_GP_bins} demonstrates the scalability of the reformulated GP \eqref{eq:GP_SSTA} on a ladder of maxima and convolutions, parameterized by the number of bins per gate.  Notice that the numbers of cones (blue line) and variables (orange line) are quadratic in the number of bins. At the same time, the relative error decreases with the number of bins, as expected. 

\vskip 6pt

\section{Discussion and Conclusions}
\label{sec:discussion}

\textcolor{black}
{In this paper, the problem of the calculation of the maximum delay in a digital circuit under uncertainty (also known as SSTA) is studied from the mathematical optimization point of view \emph{for the first time}. Using a histogram representation of the delays distributions for simplicity, we have presented two formulations of SSTA as an optimization problem. Section~\ref{sec:ssta_bip} shows Binary Integer Programming (BIP) approach, which is a formal formulation and does not scale. Section~\ref{sec:ssta_gp} gives a more practical formulation of SSTA as a Geometric Programming (GP) problem.}

For a reformulation of the GP, we have demonstrated linear scaling with the number of gates and quadratic scaling with the number of bins. The SSTA has been successfully computed using 30 bins for a circuit with 400 gates in 440 seconds \color{black} which ran on an 8-processor machine equipped with Intel Xeon Scalable Platinum 8160 (192 cores at 2.1GHz and 384 hardware threads) with 1536~GB RAM\color{black}.\footnote{Full code used in this research can be found in the repository \url{https://github.com/bosakad/GP-Optimization/}.}

\textcolor{black}{The histogram approximation, which has been previously studied~\cite{liou2001FastStatisticalTiming} as a replacement of Monte Carlo simulations, is used in this work for optimization purposes. However, this approach} has clear disadvantages: (i) as we increase the number of gates, we have to increase the size of the interval, and with that the number of bins; \textcolor{black}{(ii) we also assume the support of the delay distribution is known\footnote{\color{black} The delay in a circuit is clearly bounded from below and may well be bounded from above by considerations of practicality with respect to clocking frequency. Principled methods for establishing upper bounds have been proposed \cite[Section 2.3]{liou2001FastStatisticalTiming}, but we have used trial and error to set the upper bound in the proof-of-concept implementation in this paper.} (iii) Also, it should be noted that the correlations between the delays were not taken into account.}

On the other hand, this approach allowed us to (i) perform the robust optimization of delays' distributions, unlike in other statistical approaches, where only the statistical moments are taken into account, and (ii) perform computations in polynomial time up to any fixed precision using GP. Last but not least, the histogram formulation of the SSTA makes the results transparent and easy to understand. 

\color{black}
Some of these challenges could be addressed.  
Similar to the approximation algorithm \cite[Section 3.3]{liou2001FastStatisticalTiming} of Liou~\textit{et al.}, one could address the scalability issue at the cost of 
some error by decomposing the circuit into ``supergates'', possibly hierarchically, and  
at the further cost of estimating only the tail of the delay distribution by (i) ``filtering out unnecessary stems'', i.e., discarding sample paths, which are guaranteed not to influence the tail of the distribution.
See \cite[Section 2.3]{liou2001FastStatisticalTiming} for suggestions how this could be performed. 
The same approach of \cite[Section 2.3]{liou2001FastStatisticalTiming} could also be used to address the challenge (ii) above, as it can be used to estimate the range of arrival times of events. 
The challenge (iii) above, phrased in terms of addressing the correlations seems to be inherently difficult, albeit perhaps less important. 
This inherent difficulty extends to measuring the correlations between more than two gates' delays,
especially when conditioned on external factors such as a change of temperature.
Indeed, many sources \cite[e.g., p. 131]{jyu1993StatisticalTimingAnalysis} claim ``cell library designers agree that it is reasonable to expect the delays for components on a single chip to track each other.''

\begin{table}[t!]
\color{black}
\begin{tabular}{l|c|c}
    Model             & Simulation refs.\ & Optimization formulation \\
    \hline 
    Histogram             & Liou et al.\ \cite{liou2001FastStatisticalTiming} & Geometric program of Section \ref{sec:ssta_gp} \\
    Impulse train        & Naidu \cite{naidu2002TimingYieldCalculation}                & Mixed-integer linear \\ 
    Gaussian comb         & Mishagli et al. \cite{mishagli2020RadialBasisFunctions}  & Mixed-integer tame \\
    Gaussian mixtures     & Mishagli et al. \cite{mishagli2020RadialBasisFunctions}  & Mixed-integer tame \\ 
\end{tabular}
\caption{An overview of the approaches that approximate probability distributions with various parametric classes of distributions, with references to their uses in replacements of Monte Carlo in simulation, and our suggestions as to
the classes of optimization problems obtained using the technique of Section~\ref{sec:ssta_bip}.}
\end{table}

Let us now chart some potential avenues for further research. 
Easily, one could replace the uniform distribution centered at the midpoint of each bin
with a triangular distribution centered at the midpoint of each bin. 
In his pioneering paper \cite{naidu2002TimingYieldCalculation}, Naidu used such ``impulse train'' distributions as a replacement of Monte Carlo for simulation purposes. 
Our binary--integer optimization formulation of Section \ref{sec:ssta_bip} should be easy to extend to the triangular distributions,
and would remain mixed-integer linear.
Whether the geometric--programming approach of Section \ref{sec:ssta_gp} would be as easy to extend, remains to be investigated.
Either way, this would be an interesting extension. 

Further, our work could be plausibly extended to the Gaussian comb model of Mishagli and Blokhina \cite{mishagli2020RadialBasisFunctions}, where one would replace the uniform distribution centered at the midpoint of each bin
with a Gaussian kernel function. As such, the Gaussian comb model is a very special case
of the Gaussian mixture model with predefined, uniformly distributed expectations of the components. 
There, one optimizes over the mixture coefficients, rather than the counts in a histogram,
and the objective function has a more complicated form (which expresses the integral in a closed-form). 
We conjecture the resulting mixed-integer non-linear optimization problems are ``tame'' in their continuous
part, \textit{i.e.}, the continuous part is definable in an o-minimal structure. 
This is a fast-growing area of optimization, due to the applications in deep learning (see, \textit{e.g.},~\cite{bolte2021ConservativeSetValued}),
but mixed-integer extensions do not seem to have been studied yet, and there certainly are no off-the-shelf solvers. 

More broadly, one could consider Gaussian mixture models or infinitely-smooth radial basis functions (RBFs). 
Infinitely smooth RBF, such as Multiquadric RBF or Inverse quadratic RBF, are real-analytic ($C^{\infty }(\mathbb {R})$), 
and hence one can formulate mixed-integer analytical optimization problems over these. 
While optimization over real-analytic functions is becoming better understood \cite{absil2006StableEquilibriumPoints,kurdyka2000ProofGradientConjecture},
these optimization problems seem very challenging.
While the famous result of Kurdyka~\textit{et al.} \cite{kurdyka2000ProofGradientConjecture} shows the finite length of the gradient flows, 
local minima are not necessarily stable equilibria of the gradient-descent system, and vice versa \cite[Proposition 2]{absil2006StableEquilibriumPoints}.
That is: local minimality is neither necessary nor sufficient for stability.
Integer analytic optimization hence seems difficult to work with, other than using spatial branch-and-bound \cite{smith1999SymbolicReformulationSpatial}, which may not be sufficiently scalable.


 
 \color{black}
Future work may also involve 
an extension of the maximum computation and the gate sizing program for the case of correlated random variables and related problems of circuit design.


\section*{Declarations}

\paragraph*{Funding \;}
The research that led to these results received funding from OP RDE under Grant Agreement No CZ.02.1.01/0.0/0.0/16\_019/0000765.

\paragraph*{Conflicts of interest/Competing interests \;}
The authors have no conflicts of interest to declare that are relevant to the content of this article.


\bibliographystyle{spmpsci}      
\bibliography{ssta_lib,optimization}

\begin{thebibliography}{10}
\providecommand{\url}[1]{{#1}}
\providecommand{\urlprefix}{URL }
\expandafter\ifx\csname urlstyle\endcsname\relax
  \providecommand{\doi}[1]{DOI~\discretionary{}{}{}#1}\else
  \providecommand{\doi}{DOI~\discretionary{}{}{}\begingroup
  \urlstyle{rm}\Url}\fi

\bibitem{absil2006StableEquilibriumPoints}
Absil, P.A., Kurdyka, K.: On the stable equilibrium points of gradient systems.
\newblock Systems \& Control Letters \textbf{55}(7), 573--577 (2006).
\newblock \doi{10.1016/j.sysconle.2006.01.002}.
\newblock
  \urlprefix\url{https://www.sciencedirect.com/science/article/pii/S0167691106000107}

\bibitem{agrawal2018RewritingSystemConvex}
Agrawal, A., Verschueren, R., Diamond, S., Boyd, S.: A rewriting system for
  convex optimization problems.
\newblock Journal of Control and Decision \textbf{5}(1), 42--60 (2018).
\newblock \doi{10.1080/23307706.2017.1397554}.
\newblock \urlprefix\url{https://doi.org/10.1080/23307706.2017.1397554}.
\newblock Publisher: Taylor \& Francis \_eprint:
  https://doi.org/10.1080/23307706.2017.1397554

\bibitem{beece2014TransistorSizingCustom}
Beece, D.K., Visweswariah, C., Xiong, J., Zolotov, V.: Transistor sizing of
  custom high-performance digital circuits with parametric yield
  considerations.
\newblock Optimization and Engineering \textbf{15}(1), 217--241 (2014).
\newblock \doi{10.1007/s11081-012-9208-0}.
\newblock \urlprefix\url{https://doi.org/10.1007/s11081-012-9208-0}

\bibitem{berkelaar1990GateSizingMOS}
Berkelaar, M.R.C.M., Jess, J.A.G.: Gate sizing in {MOS} digital circuits with
  linear programming.
\newblock In: Proceedings of the conference on {European} design automation,
  {EURO}-{DAC} '90, pp. 217--221. IEEE Computer Society Press, Washington, DC,
  USA (1990)

\bibitem{bhasker2009StaticTimingAnalysis}
Bhasker, J., Chadha, R.: Static {Timing} {Analysis} for {Nanometer} {Designs}.
  {A} {Practical} {Approach}.
\newblock Springer (2009)

\bibitem{blaauw2008StatisticalTimingAnalysis}
Blaauw, D., Chopra, K., Srivastava, A., Scheffer, L.: Statistical timing
  analysis: {From} basic principles to state of the art.
\newblock IEEE Trans. Comput.–Aided Des. Integr. Circuits Syst. \textbf{4}(8)
  (2008).
\newblock \doi{10.1109/TCAD.2007.907047}

\bibitem{bolte2021ConservativeSetValued}
Bolte, J., Pauwels, E.: Conservative set valued fields, automatic
  differentiation, stochastic gradient methods and deep learning.
\newblock Mathematical Programming \textbf{188}(1), 19--51 (2021).
\newblock \doi{10.1007/s10107-020-01501-5}.
\newblock \urlprefix\url{https://doi.org/10.1007/s10107-020-01501-5}

\bibitem{boyd2007TutorialGeometricProgramming}
Boyd, S., Kim, S.J., Vandenberghe, L., Hassibi, A.: A tutorial on geometric
  programming.
\newblock Optimization and Engineering \textbf{8}(1), 67--127 (2007).
\newblock \doi{10.1007/s11081-007-9001-7}.
\newblock \urlprefix\url{http://link.springer.com/10.1007/s11081-007-9001-7}

\bibitem{boyd2005DigitalCircuitOptimization}
Boyd, S.P., Kim, S.J., Patil, D.D., Horowitz, M.A.: Digital {Circuit}
  {Optimization} via {Geometric} {Programming}.
\newblock Operations Research \textbf{53}(6), 899--932 (2005).
\newblock \doi{10.1287/opre.1050.0254}.
\newblock
  \urlprefix\url{https://pubsonline.informs.org/doi/abs/10.1287/opre.1050.0254}.
\newblock Publisher: INFORMS

\bibitem{brenner2009AnalyticalMethodsPlacement}
Brenner, U., Vygen, J.: Analytical {Methods} in {Placement}.
\newblock In: Handbook of {Algorithms} for {Physical} {Design} {Automation}.
  Auerbach Publications (2009).
\newblock Num Pages: 20

\bibitem{chang2005StatisticalTimingAnalysis}
Chang, H., Sapatnekar, S.: Statistical timing analysis under spatial
  correlations.
\newblock IEEE Transactions on Computer-Aided Design of Integrated Circuits and
  Systems \textbf{24}(9), 1467--1482 (2005).
\newblock \doi{10.1109/TCAD.2005.850834}.
\newblock Conference Name: IEEE Transactions on Computer-Aided Design of
  Integrated Circuits and Systems

\bibitem{chang2003StatisticalTimingAnalysis}
Chang, H., Sapatnekar, S.S.: Statistical timing analysis considering spatial
  correlations using a single {PERT}-like traversal.
\newblock In: Proc. {ICCAD}, pp. 621--625 (2003).
\newblock \doi{10.1109/ICCAD.2003.159746}

\bibitem{chang2005ParameterizedBlockbasedStatistical}
Chang, H., Zolotov, V., Narayan, S., Visweswariah, C.: Parameterized
  block-based statistical timing analysis with non-gaussian parameters,
  nonlinear delay functions.
\newblock In: Proceedings of the 42nd annual {Design} {Automation}
  {Conference}, {DAC} '05, pp. 71--76. Association for Computing Machinery, New
  York, NY, USA (2005).
\newblock \doi{10.1145/1065579.1065604}.
\newblock \urlprefix\url{https://doi.org/10.1145/1065579.1065604}

\bibitem{cheng2012FourierSeriesApproximation}
Cheng, L., Gong, F., Xu, W., Xiong, J., He, L., Sarrafzadeh, M.: Fourier
  {Series} {Approximation} for {Max} {Operation} in {Non}-{Gaussian} and
  {Quadratic} {Statistical} {Static} {Timing} {Analysis}.
\newblock IEEE Transactions on Very Large Scale Integration (VLSI) Systems
  \textbf{20}(8), 1383--1391 (2012).
\newblock \doi{10.1109/TVLSI.2011.2157843}

\bibitem{cheng2007NonLinearStatisticalStatic}
Cheng, L., Xiong, J., He, L.: Non-{Linear} {Statistical} {Static} {Timing}
  {Analysis} for {Non}-{Gaussian} {Variation} {Sources}.
\newblock In: 2007 44th {ACM}/{IEEE} {Design} {Automation} {Conference}, pp.
  250--255 (2007).
\newblock ISSN: 0738-100X

\bibitem{clark1961GreatestFiniteSet}
Clark, C.E.: The {Greatest} of a {Finite} {Set} of {Random} {Variables}.
\newblock Oper. Res. \textbf{9}(2), 145--162 (1961).
\newblock \doi{10.1287/opre.9.2.145}.
\newblock Place: Institute for Operations Research and the Management Sciences
  (INFORMS), Linthicum, Maryland, USA Publisher: INFORMS

\bibitem{diamond2016CVXPYPythonembeddedModeling}
Diamond, S., Boyd, S.: {CVXPY}: a python-embedded modeling language for convex
  optimization.
\newblock The Journal of Machine Learning Research \textbf{17}(1), 2909--2913
  (2016)

\bibitem{forzan2009StatisticalStaticTiming}
Forzan, C., Pandini, D.: Statistical static timing analysis: {A} survey.
\newblock Integration, the VLSI Journal \textbf{42}(3), 409--435 (2009).
\newblock \doi{10.1016/j.vlsi.2008.10.002}

\bibitem{freeley2018StatisticalSimulationsDelay}
Freeley, J., Mishagli, D., Brazil, T., Blokhina, E.: Statistical {Simulations}
  of {Delay} {Propagation} in {Large} {Scale} {Circuits} {Using} {Graph}
  {Traversal} and {Kernel} {Function} {Decomposition}.
\newblock In: Proc. {SMACD} (2018)

\bibitem{held2011CombinatorialOptimizationVLSI}
Held, S., Korte, B., Rautenbach, D., Vygen, J.: Combinatorial {Optimization} in
  {VLSI} {Design}.
\newblock Combinatorial Optimization pp. 33--96 (2011).
\newblock \doi{10.3233/978-1-60750-718-5-33}.
\newblock \urlprefix\url{http://www.or.uni-bonn.de/research/montreal.pdf}.
\newblock Publisher: IOS Press

\bibitem{jacobs2000GateSizingUsing}
Jacobs, E.T.A.F., Berkelaar, M.R.C.M.: Gate sizing using a statistical delay
  model.
\newblock In: Proceedings {Design}, {Automation} and {Test} in {Europe}
  {Conference} and {Exhibition} 2000, pp. 283--290 (2000).
\newblock \doi{10.1109/DATE.2000.840285}

\bibitem{jess2003StatisticalTimingParametric}
Jess, J.A.G., Kalafala, K., Naidu, S.R., Otten, R.H.J.M., Visweswariah, C.:
  Statistical {Timing} for {Parametric} {Yield} {Prediction} of {Digital}
  {Integrated} {Circuits}.
\newblock In: Proc. {DAC}, pp. 932--937. ACM (2003).
\newblock \doi{10.1145/775832.776066}.
\newblock Event-place: Anaheim, CA, USA

\bibitem{jess2006StatisticalTimingParametric}
Jess, J.A.G., Kalafala, K., Naidu, S.R., Otten, R.H.J.M., Visweswariah, C.:
  Statistical {Timing} for {Parametric} {Yield} {Prediction} of {Digital}
  {Integrated} {Circuits}.
\newblock IEEE Trans. Comput.–Aided Des. Integr. Circuits Syst.
  \textbf{25}(11), 2376--2392 (2006).
\newblock \doi{10.1109/TCAD.2006.881332}

\bibitem{joshi2008EfficientMethodLargeScale}
Joshi, S., Boyd, S.: An {Efficient} {Method} for {Large}-{Scale} {Gate}
  {Sizing}.
\newblock IEEE Transactions on Circuits and Systems I: Regular Papers
  \textbf{55}(9), 2760--2773 (2008).
\newblock \doi{10.1109/TCSI.2008.920087}.
\newblock Conference Name: IEEE Transactions on Circuits and Systems I: Regular
  Papers

\bibitem{jyu1993StatisticalTimingAnalysis}
Jyu, H.F., Malik, S., Devadas, S., Keutzer, K.: Statistical timing analysis of
  combinational logic circuits.
\newblock IEEE Transactions on Very Large Scale Integration (VLSI) Systems
  \textbf{1}(2), 126--137 (1993).
\newblock \doi{10.1109/92.238423}.
\newblock Conference Name: IEEE Transactions on Very Large Scale Integration
  (VLSI) Systems

\bibitem{khandelwal2005GeneralFrameworkAccurate}
Khandelwal, V., Srivastava, A.: A general framework for accurate statistical
  timing analysis considering correlations.
\newblock In: Proc. {DAC}, pp. 89--94 (2005).
\newblock \doi{10.1145/1065579.1065607}.
\newblock ISSN: 0738-100X

\bibitem{knuth1993JohannFaulhaberSums}
Knuth, D.E.: Johann {Faulhaber} and sums of powers.
\newblock Mathematics of Computation \textbf{61}(203), 277--294 (1993).
\newblock \doi{10.1090/S0025-5718-1993-1197512-7}.
\newblock
  \urlprefix\url{https://www.ams.org/mcom/1993-61-203/S0025-5718-1993-1197512-7/}

\bibitem{korte2008CombinatorialProblemsChip}
Korte, B., Vygen, J.: Combinatorial {Problems} in {Chip} {Design}.
\newblock In: M.~Grötschel, G.O.H. Katona, G.~Sági (eds.) Building {Bridges}:
  {Between} {Mathematics} and {Computer} {Science}, Bolyai {Society}
  {Mathematical} {Studies}, pp. 333--368. Springer, Berlin, Heidelberg (2008).
\newblock \doi{10.1007/978-3-540-85221-6_12}.
\newblock \urlprefix\url{https://doi.org/10.1007/978-3-540-85221-6_12}

\bibitem{kurdyka2000ProofGradientConjecture}
Kurdyka, K., Mostowski, T., Parusinski, A.: Proof of the {Gradient}
  {Conjecture} of {R}. {Thom}.
\newblock Annals of Mathematics \textbf{152}(3), 763--792 (2000).
\newblock \doi{10.2307/2661354}.
\newblock \urlprefix\url{https://www.jstor.org/stable/2661354}.
\newblock Publisher: Annals of Mathematics

\bibitem{liou2001FastStatisticalTiming}
Liou, J.J., Cheng, K.T., Kundu, S., Krstic, A.: Fast statistical timing
  analysis by probabilistic event propagation.
\newblock In: Proceedings of the 38th {Design} {Automation} {Conference}
  ({IEEE} {Cat}. {No}.{01CH37232}), pp. 661--666 (2001).
\newblock \doi{10.1145/378239.379043}.
\newblock ISSN: 0738-100X

\bibitem{mishagli2020RadialBasisFunctions}
Mishagli, D., Blokhina, E.: Radial {Basis} {Functions} {Based} {Algorithms} for
  {Non}-{Gaussian} {Delay} {Propagation} in {Very} {Large} {Circuits}.
\newblock In: V.V. Krzhizhanovskaya, G.~Závodszky, M.H. Lees, J.J. Dongarra,
  P.M.A. Sloot, S.~Brissos, J.~Teixeira (eds.) Computational {Science} –
  {ICCS} 2020, Lecture {Notes} in {Computer} {Science}, pp. 217--229. Springer
  International Publishing, Cham (2020).
\newblock \doi{10.1007/978-3-030-50426-7_17}

\bibitem{naidu2002TimingYieldCalculation}
Naidu, S.: Timing yield calculation using an impulse-train approach.
\newblock In: Proceedings of {ASP}-{DAC}/{VLSI} {Design} 2002. 7th {Asia} and
  {South} {Pacific} {Design} {Automation} {Conference} and 15h {International}
  {Conference} on {VLSI} {Design}, pp. 219--224 (2002).
\newblock \doi{10.1109/ASPDAC.2002.994923}

\bibitem{naidu2021ConvexProgrammingSolution}
Naidu, S.R.: A convex programming solution for gate-sizing with pipelining
  constraints.
\newblock Optimization and Engineering  (2021).
\newblock \doi{10.1007/s11081-021-09616-0}.
\newblock \urlprefix\url{https://doi.org/10.1007/s11081-021-09616-0}

\bibitem{rakai2015SizingDigitalCircuits}
Rakai, L., Farshidi, A.: Sizing {Digital} {Circuits} {Using} {Convex}
  {Optimization} {Techniques}.
\newblock In: M.~Fakhfakh, E.~Tlelo-Cuautle, P.~Siarry (eds.) Computational
  {Intelligence} in {Digital} and {Network} {Designs} and {Applications}, pp.
  3--32. Springer International Publishing, Cham (2015).
\newblock \doi{10.1007/978-3-319-20071-2_1}.
\newblock \urlprefix\url{https://doi.org/10.1007/978-3-319-20071-2_1}

\bibitem{ramprasath2016SkewNormalCanonicalModel}
Ramprasath, S., Vijaykumar, M., Vasudevan, V.: A {Skew}-{Normal} {Canonical}
  {Model} for {Statistical} {Static} {Timing} {Analysis}.
\newblock IEEE Transactions on Very Large Scale Integration (VLSI) Systems
  \textbf{24}(6), 2359--2368 (2016).
\newblock \doi{10.1109/TVLSI.2015.2501370}

\bibitem{rockafellar2000optimization}
Rockafellar, R.T., Uryasev, S., et~al.: Optimization of conditional
  value-at-risk.
\newblock Journal of risk \textbf{2}, 21--42 (2000)

\bibitem{sapatnekar2004Timing}
Sapatnekar, S.: Timing.
\newblock Springer-Verlag (2004)

\bibitem{sapatnekar1993ExactSolutionTransistor}
Sapatnekar, S., Rao, V., Vaidya, P., Kang, S.M.: An exact solution to the
  transistor sizing problem for {CMOS} circuits using convex optimization.
\newblock IEEE Transactions on Computer-Aided Design of Integrated Circuits and
  Systems \textbf{12}(11), 1621--1634 (1993).
\newblock \doi{10.1109/43.248073}.
\newblock Conference Name: IEEE Transactions on Computer-Aided Design of
  Integrated Circuits and Systems

\bibitem{singh2008ScalableStatisticalStatic}
Singh, J., Sapatnekar, S.S.: A {Scalable} {Statistical} {Static} {Timing}
  {Analyzer} {Incorporating} {Correlated} {Non}-{Gaussian} and {Gaussian}
  {Parameter} {Variations}.
\newblock IEEE Trans. Comput.–Aided Des. Integr. Circuits Syst.
  \textbf{27}(1), 160--173 (2008).
\newblock \doi{10.1109/TCAD.2007.907241}

\bibitem{smith1999SymbolicReformulationSpatial}
Smith, E.M.B., Pantelides, C.C.: A symbolic reformulation/spatial
  branch-and-bound algorithm for the global optimisation of nonconvex {MINLPs}.
\newblock Computers \& Chemical Engineering \textbf{23}(4), 457--478 (1999).
\newblock \doi{10.1016/S0098-1354(98)00286-5}.
\newblock
  \urlprefix\url{https://www.sciencedirect.com/science/article/pii/S0098135498002865}

\bibitem{vielma2011modeling}
Vielma, J.P., Nemhauser, G.L.: Modeling disjunctive constraints with a
  logarithmic number of binary variables and constraints.
\newblock Mathematical Programming \textbf{128}, 49--72 (2011)

\bibitem{visweswariah2003DeathTaxesFailing}
Visweswariah, C.: Death, taxes and failing chips.
\newblock In: Proc. {DAC}, pp. 343--347. IEEE (2003).
\newblock \doi{10.1145/775832.775921}.
\newblock Event-place: Anaheim, CA, USA

\bibitem{visweswariah2004FirstOrderIncrementalBlockBased}
Visweswariah, C., Ravindran, K., Kalafala, K., Walker, S.G., Narayan, S.:
  First-{Order} {Incremental} {Block}-{Based} {Statistical} {Timing}
  {Analysis}.
\newblock In: Proc. {DAC}, pp. 331--336. ACM, New York, NY, USA (2004).
\newblock \doi{10.1145/996566.996663}.
\newblock Event-place: San Diego, CA, USA

\bibitem{visweswariah2006FirstOrderIncrementalBlockBased}
Visweswariah, C., Ravindran, K., Kalafala, K., Walker, S.G., Narayan, S.,
  Beece, D.K., Piaget, J., Venkateswaran, N., Hemmett, J.G.: First-{Order}
  {Incremental} {Block}-{Based} {Statistical} {Timing} {Analysis}.
\newblock IEEE Trans. Comput.–Aided Des. Integr. Circuits Syst.
  \textbf{25}(10), 2170--2180 (2006).
\newblock \doi{10.1109/TCAD.2005.862751}

\bibitem{vygen2006SlackStaticTiming}
Vygen, J.: Slack in static timing analysis.
\newblock IEEE Transactions on Computer-Aided Design of Integrated Circuits and
  Systems \textbf{25}(9), 1876--1885 (2006).
\newblock \doi{10.1109/TCAD.2005.858348}.
\newblock Conference Name: IEEE Transactions on Computer-Aided Design of
  Integrated Circuits and Systems

\bibitem{wolsey1999IntegerCombinatorialOptimization}
Wolsey, L., Nemhauser, G.: Integer and combinatorial optimization.
\newblock Wiley series in discrete mathematics and optimization. Wiley (1999).
\newblock
  \urlprefix\url{https://www.wiley.com/en-us/Integer+and+Combinatorial+Optimization-p-9780471359432}

\bibitem{zhan2005CorrelationawareStatisticalTiming}
Zhan, Y., Strojwas, A.J., Li, X., Pileggi, L.T., Newmark, D., Sharma, M.:
  Correlation-aware statistical timing analysis with non-{Gaussian} delay
  distributions.
\newblock In: Proc. {DAC}, pp. 77--82 (2005).
\newblock \doi{10.1145/1065579.1065605}.
\newblock ISSN: 0738-100X

\bibitem{zhang2005CorrelationpreservedNonGaussianStatistical}
Zhang, L., Chen, W., Hu, Y., Gubner, J.A., Chen, C.C.P.: Correlation-preserved
  non-{Gaussian} statistical timing analysis with quadratic timing model.
\newblock In: Proc. {DAC}, pp. 83--88 (2005).
\newblock \doi{10.1109/DAC.2005.193778}.
\newblock ISSN: 0738-100X

\bibitem{zhang2006StatisticalTimingAnalysis}
Zhang, L., Hu, Y., Chen, C.C.P.: Statistical timing analysis with path
  reconvergence and spatial correlations.
\newblock In: Proc. {DAC}, vol.~1, p. 5 pp. (2006).
\newblock \doi{10.1109/DATE.2006.243890}.
\newblock ISSN: 1530-1591

\end{thebibliography}

\appendix
\FloatBarrier

\section{Note on unary encoding}

The idea behind the unary encoding is as follows. A histogram, which is an array of real numbers, can be written as a matrix of \emph{binary numbers}. Each bin is then expressed by a row in the matrix. In this case, the probability is given by a sum of the row elements divided by both a sum of the matrix elements and a width of the bin (for normalization purposes).

\begin{figure}[t]
\centering
\includegraphics[scale=0.6]{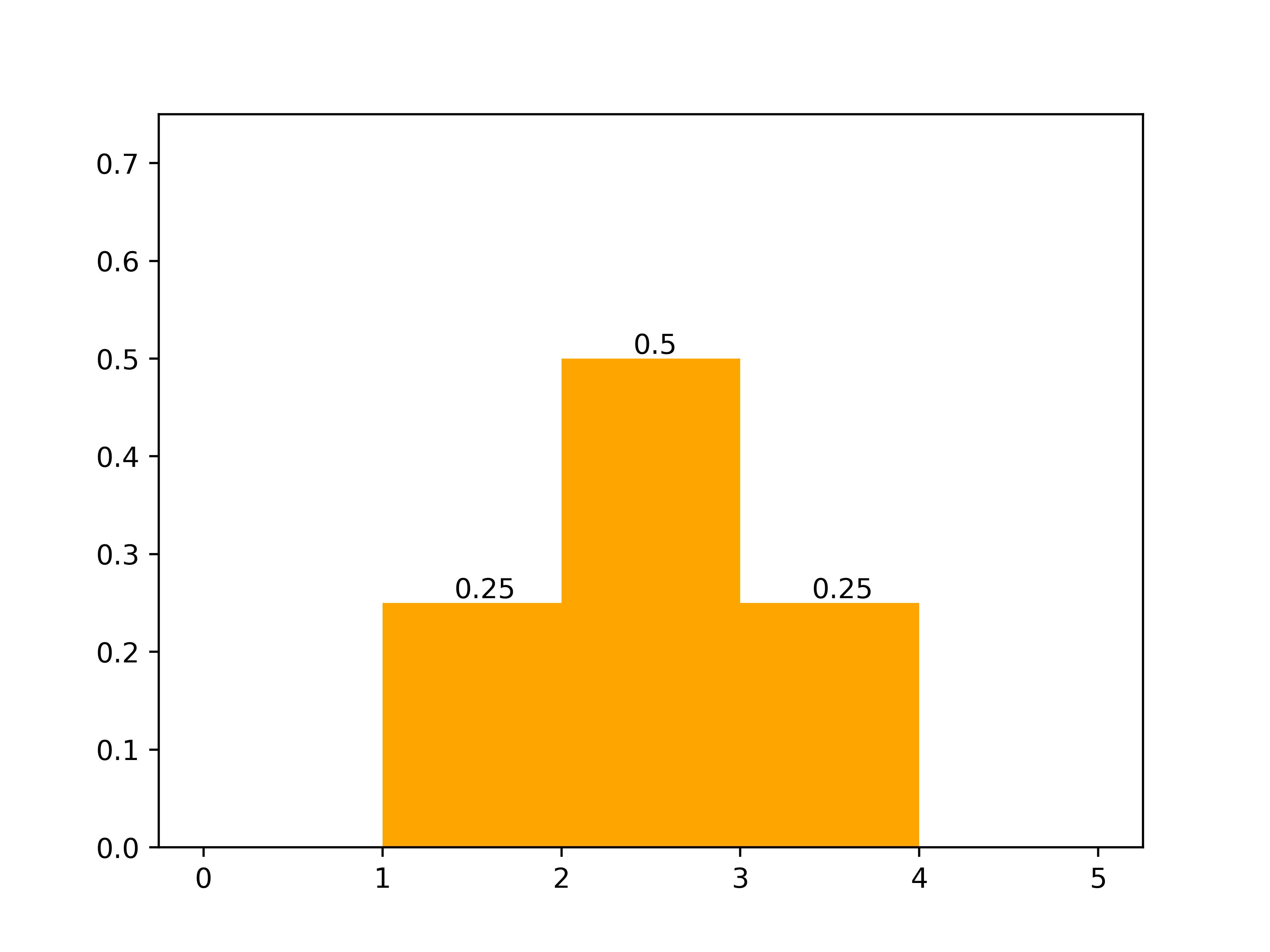}
\caption{An example histogram with $n=5$ bins; width of the bin is 1.}
\label{histogramUnary}
\end{figure}

Consider a toy example shown in Fig.~\ref{histogramUnary}. This histogram can be represented by a $5\times3$ matrix as follows:

\begin{equation*}\label{matrixUnary}
\mathbf H = 
\begin{bmatrix}
0 & 0 & 0\\
1 & 0 & 0\\
1 & 1 & 0\\
1 & 0 & 0\\
0 & 0 & 0
\end{bmatrix}.
\end{equation*}
Note that not all real numbers can be encoded by a finite number of binary numbers. The accuracy of the encoding is proportional to the number of columns in the binary representation.

{\color{black}
\paragraph{Numerical example.} Consider two histograms:
\begin{equation}\label{eq:histograms_example}
\mathbf H_A = 
\begin{bmatrix}
0 & 0 & 0\\
1 & 0 & 0\\
1 & 1 & 0\\
1 & 0 & 0\\
0 & 0 & 0
\end{bmatrix}
\quad
\mathbf H_B = 
\begin{bmatrix}
0 & 0 & 0\\
0 & 0 & 0\\
1 & 1 & 1\\
1 & 0 & 0\\
0 & 0 & 0
\end{bmatrix}
\end{equation}

The convolution $\mathbf H_{AB}$ of these histograms, according to the Algorithm\ref{alg:convolution_bip}, will require introducing the constraints \eqref{binary_constr} for the unary numbers in $\mathbf H_A$ and $\mathbf H_B$. For example, for $z=3$ we have $k=1,2$, and there will be the following values
\begin{align*}
\mathbf H_A[k,i] &= \{ \mathbf H_A[1,i], \mathbf H_A[2,i] \}, & i &= 1,2,3
\\
\mathbf H_B[z-k,i] &= \{ \mathbf H_B[2,j], \mathbf H_B[1,j] \}, & j &= 1,2,3
\end{align*}
The constraints \eqref{binary_constr} are needed to be introduced for the permutation of these sets. Then, the constraints are passed to a solver together with the auxiliary variables $s$.
}

\end{document}